\documentclass[]{imsart}

\usepackage{amsthm,amsmath}
\usepackage{amssymb}
\usepackage{amsfonts}
\usepackage{multirow,dsfont,mathtools}
\RequirePackage[numbers]{natbib}
\RequirePackage[colorlinks,citecolor=blue,urlcolor=blue]{hyperref}
\usepackage{siunitx}
\usepackage{graphicx}
\usepackage{subcaption}

\startlocaldefs

\newtheorem{theorem}{Theorem}
\numberwithin{theorem}{section}
\newtheorem{lemma}{Lemma}
\numberwithin{lemma}{section}
\newtheorem{corollary}{Corollary}
\numberwithin{corollary}{section}
\newtheorem{proposition}{Proposition}
\numberwithin{proposition}{section}
\theoremstyle{definition}

\numberwithin{definition}{section}

\numberwithin{example}{section}
\newtheorem{remark}{Remark}
\numberwithin{remark}{section}

\newcommand{\cov}{\mathbf{Cov}}

\newcommand{\R}{\mathbb{R}}
\newcommand{\E}{\mathbb{E}}
\newcommand{\var}{\mathbf{Var}}
\newcommand{\N}{\mathbb{N}}

\newcommand{\vc}{\mathbf{1}}

\renewcommand{\P}{\mathbb{P}}

\endlocaldefs

\begin{document}

\begin{frontmatter}

\title{Statistical inference on $D^{(d)}(u_n)$ condition and estimation of the Extremal Index}
\runtitle{On $D^{(d)}(u_n)$ condition and the extremal index}
\runauthor{J.J. Cai}
\begin{aug}
\author{\fnms{Juan-Juan} \snm{Cai}\ead[label=e1]{j.cai@vu.nl}}
\address{Department of Econometrics and Data Science, Vrije Universiteit Amsterdam\\ De Boelelaan 1105, 1081HV, Amsterdam, the Netherlands\\
\printead{e1}}
\end{aug}


\begin{abstract}
Clustering of extreme events can have profound and detrimental societal consequences. The extremal index, a number in the unit interval, is a key parameter in modelling the clustering of extremes. The study of extremal index often assumes a local dependence condition known as the $D^{(d)}(u_n)$ condition.
In this paper, we develop a hypothesis test for $D^{(d)}(u_n)$ condition based on asymptotic results. We develop an estimator for the extremal index by leveraging the inference procedure based on the $D^{(d)}(u_n)$ condition, and we establish the asymptotic normality of this estimator. The finite sample performances of the hypothesis test and the estimation are examined in a simulation study, where we  consider both models that satisfies the $D^{(d)}(u_n)$ condition and models that violate this condition. In a simple case study, our statistical procedure shows that  daily temperature in summer shares  a common clustering structure of extreme values based on the data observed in three weather stations in the Netherlands, Belgium and Spain.


\end{abstract}





\end{frontmatter}

\section{Introduction}\label{Sec:Introduction}

\begin{figure}[t!]
\begin{center}
\includegraphics[width=0.85\textwidth]{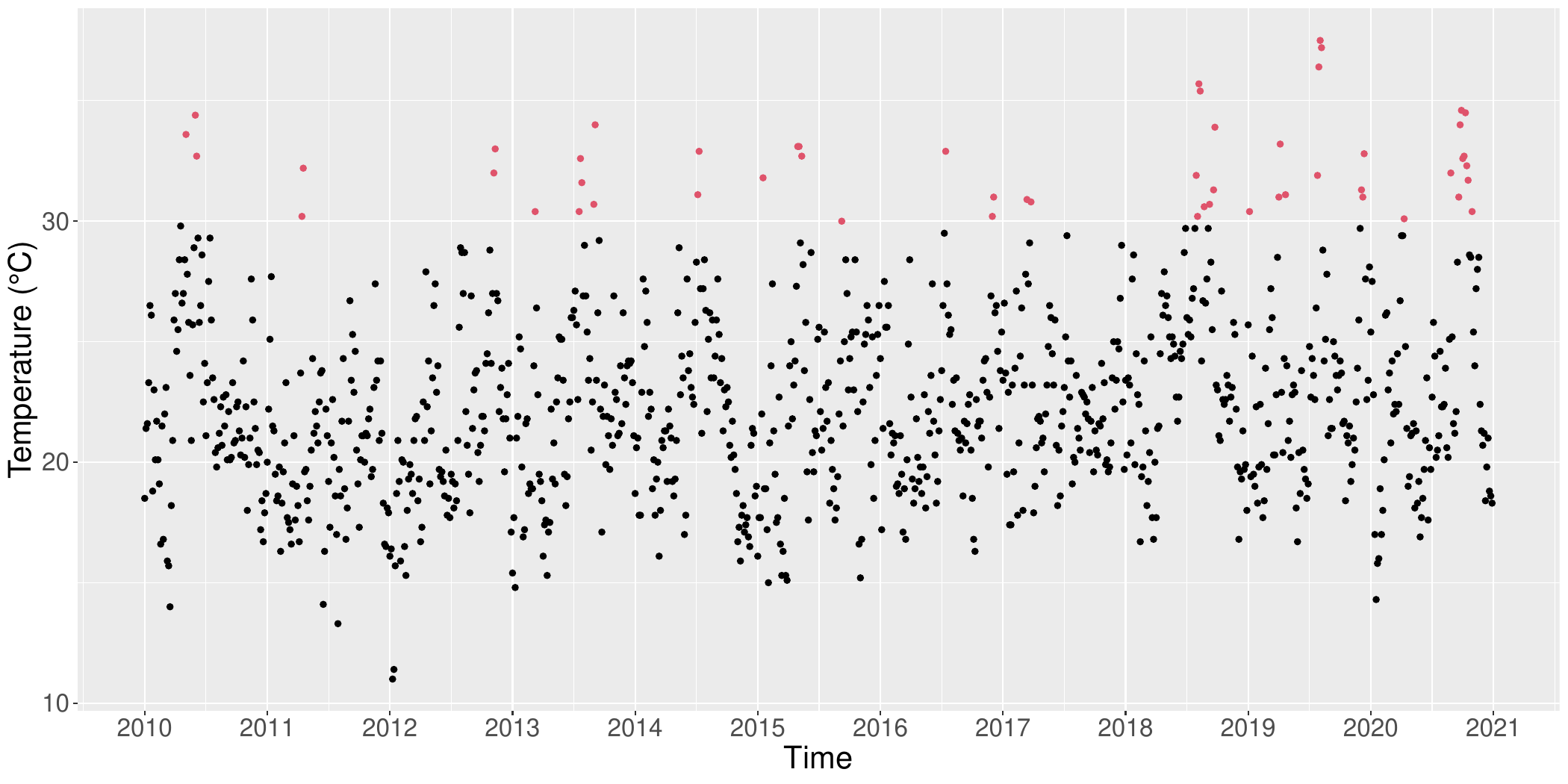}
\caption{Daily maximum temperature in June, July and
August, measured in de Bilt, the Netherlands}
\label{fig: tem_peaks}
\end{center}
\end{figure}

A cluster of extremes refers to the occurrence of multiple extreme observations (i.e.\ high
level exceedances) within a short period of time. When extreme events happen sequentially, it often has a destructive impact on our society. For instance, hot temperature extremes in successive days increase the risk of mortality, drought, wildfire and others.  Figure \ref{fig: tem_peaks} shows a typical clustering behaviour of hot days in summer in the Netherlands. The clustering of extremes is due to the serial extremal dependence of time series data. The extremal index $\theta$  is a key parameter to measure the strength of such dependence. 
The smaller value of $\theta$ indicates stronger extremal dependence, while $\theta=1$ corresponds to the case where the extremes are (asymptotically) independent, meaning that there is no clustering of extremes. The extremal index  provides two other insights on the behavior of the clustering. First, it equals the reciprocal of the expected clustering size,  that is the number of extreme observations in a cluster (\cite{Leadbetter1983}).
	Second, \cite{OBrien1987} shows that $\theta$ equals a conditional probability that measures to what extent extremes cluster together. We refer to \cite{Moloneyetal2019} for an overview on the probability properties and the accompanying interpretations of the extremal index. The primary goal of this paper is to develop an estimation for the extremal index. 

The existing methods for estimating $\theta$ can be grouped into two categories. The first category of estimators requires two tuning parameters, namely a threshold and a cluster size parameter. The threshold indicates a level above which observations are considered as extremes.  The cluster size parameter specifies the time difference between two extreme events that can be treated as independent.  Both parameters depend on the underlying asymptotic properties of the time series and thus are difficult to choose in practice.  The representatives of this group are the blocks and runs estimators (\cite{SmithWeissman1994, WeissmanNovak1998}),  which are based on the two aforementioned interpretations of $\theta$, respectively.    The second category of estimation methods need only one tuning parameter.  The inter-exceedance times method developed by \cite{FerroSegers2003} requires only the threshold parameter. The maximum likelihood estimators (MLE) of $\theta$  studied by  \cite{Northrop2015} and \cite{BerghausBucher2018} needs the cluster size parameter. Note that this MLE is a nonparametric method and is based on the following intuition: for $(X_i)$  a stationary sequence and $F$  the marginal distribution function,  $r_n(1-F(\max_{i=1,\ldots,r_n}X_i))$  is approximately exponentially distributed with parameter $\theta$, where $r_n$ denotes the cluster size parameter. There are many other references on the estimation of extremal index such as \cite{Hsing1993}, \cite{LauriniTawn2003}, \cite{Robert2009}, and \cite{Anconaetal2000}.

In this paper, we develop an estimator of $\theta$ assuming a local dependence condition, namely the $D^{(d)}(u_n)$ condition (cf.\ Section \ref{sec: Preliminary results}), where $d$ is a positive integer and $u_n$ is the threshold parameter.  Under the $D^{(d)}(u_n)$ condition, the existence and the value of $\theta$  are fully determined by the tail dependence of $(X_1,\ldots, X_d)$. For instance, $D^{(1)}(u_n)$ implies that $\theta=1$, that is the extremes tend to occur independently. Similar variants of $D^{(1)}(u_n)$ condition were assumed in \cite{Watson1954, Loynes1965} and Chapter 3.4 in \cite{LeadbetterLindgrenRootzen1983} for stationary sequences of which the sample maximum has the same limit as that of an independent sample. A slightly stronger version of $D^{(2)}(u_n)$ was introduced in \cite{LeadbetterNandagopalan1989} to allow for clustering of extremes in a stationary sequence. \cite{Chernicketal1991} introduced the $D^{(d)}(u_n)$ condition for a general integer $d$ and stated the necessary and sufficient condition for the existence of $\theta$. Building upon these probabilistic results, several authors have developed different estimators of $\theta$ under the $D^{(d)}(u_n)$ condition. \cite{Hsing1993} has proposed an estimator of $\theta$ assuming that $d$ is known and established the asymptotic normality of the estimator under $m$-dependence. \cite{Suveges2007} has studied a likelihood estimator of $\theta$ under $D^{(2)}(u_n)$. \cite{FerreiraFerreira2018} has proposed an estimator of $\theta$ by relating a stationary sequence satisfying $D^{(d)}(u_n)$ condition to a regenerative process satisfying $D^{(1)}(u_n)$ or $D^{(2)}(u_n)$ condition. \cite{HolevsovskyFusek2020} has considered an estimator of $\theta$ by introducing a censoring to the inter-exceedance times under $D^{(d)}(u_n)$ condition. In this paper, we develop a hypothesis test for $D^{(d)}(u_n)$ condition. To the best of our knowledge, this is the first attempt to propose an inference procedure for verifying $D^{(d)}(u_n)$ condition based on asymptotic results. Next, we construct a consistent estimator of the smallest $d$ such that $D^{(d)}(u_n)$ condition is valid and use this in the estimation of $\theta$.  We prove the asymptotic normality of our proposed estimator of $\theta$, which requires only one tuning parameter, the threshold.   Among the existing literature on estimating $\theta$, only a few  have addressed the asymptotic normality, cf.\ \cite{Hsing1993, WeissmanNovak1998} and \cite{BerghausBucher2018}. Yet, our result is proved under a rather general setting compared to the existing ones. Another contribution of this paper is that we express $\theta$ using the stable tail dependence function, a commonly used tool from multivariate extreme value theory.  This result links the extremal behavior of a stationary sequence to multivariate extreme value theory and it provides a convenient way to compute the value of $\theta$ for some time series models.

 The rest of this paper is organized as follows.  The hypothesis test on the $D^{(d)}(u_n)$ condition and the estimator of $\theta$  are developed in Section \ref{sec:AP} alongside the corresponding  asymptotic properties. In Section \ref{sec:simulation}, we demonstrate the computation of $\theta$ via the stable tail dependence function for some examples, and investigate the finite sample performance of the developed hypothesis test and the estimation of $\theta$.  In Section \ref{Sec:RealData}, we apply our estimations to summer daily temperature observed in three weather stations located respectively in the Netherlands, Belgium and Spain. Our estimates show that the extremal dependence of hot days behaves very closely in these three stations despite of the different weather types. Sections \ref{Sec:Proof} and \ref{sec: lem} contain the proofs for the main theorems.

 \section{Estimations and asymptotic properties} \label{sec:AP}

\subsection{Preliminary results}  \label{sec: Preliminary results}
Let $\{X_i\}$ be strictly stationary sequence of random variables with a continuous marginal distribution function $F$. Letting $u_n(\tau)$, $\tau>0$ be such that 
\begin{equation}
\lim_{n\rightarrow\infty}n(1-F(u_n(\tau))=\tau, \label{eq: un}
\end{equation}
we say that $\{X_i\}$  has an extremal index $\theta \in [0, 1]$, if 
$$
\lim_{n\rightarrow\infty}\P\left(\max_{1\leq i\leq n}X_i\leq u_n(\tau)\right)=e^{-\theta\tau}.
$$

 Our goal is to estimate the extremal index $\theta$. The starting point is  Corollary 1.3 from \cite{Chernicketal1991}.  For stating the existence of $\theta$, the corollary assumes two  ``mixing" conditions in the dependence of the sequence, namely $D(u_n(\tau))$ and $D^{(d)}(u_n(\tau))$, where $u_n(\tau)$ is given in \eqref{eq: un}.
When it does not cause any misunderstanding, we write $u_n$ instead of  $u_n(\tau)$.

\emph{Condition $D(u_n)$} For any integers $1\leq i_1< \cdots < i_q <j_1<\cdots<j_{q'}\leq n$ for which $j_1-i_q\geq l$, we have
\begin{align*}
\left|\P\left(\max_{1\leq t\leq q}X_{i_t}\leq u_n,\max_{1\leq t\leq q'}X_{j_t}\leq u_n\right)
 -\P\left( \max_{1\leq t\leq q}X_{i_t}\leq u_n\right)\P\left(\max_{1\leq t\leq q'}X_{j_t}\leq u_n\right)\right|
\leq \alpha_{n,l},
\end{align*}
and $\lim_{n\rightarrow \infty}\alpha_{n,l_n}=0$ for some sequence $l_n=o(n)$ and $l_n \rightarrow \infty$.

\emph{Condition $D^{(d)}(u_n)$}  For a positive integer $d$,  there exist sequences of integers ${r_n}$ and $l_n$ such that $r_n\rightarrow \infty$, $n\alpha_{n,l_n}/r_n\rightarrow 0$, $l_n/r_n\rightarrow 0$, and 
\begin{equation}
\lim_{n\rightarrow \infty} n\P(X_1>u_n\geq M_{2,d}, M_{d+1,r_n}>u_n)=0,   \label{eq: DdUn}
\end{equation}
where $M_{i,j}:=-\infty$ for $i>j$ and $M_{i,j}:=\max_{i\leq t \leq j} X_t$ for $i\leq j$.

Condition $D(u_n)$ is a standard condition on long range dependence when studying the extremal behaviour of a stationary sequence (see e.g. \cite{LeadbetterLindgrenRootzen1983}). Condition $D^{(d)}(u_n)$ is a local mixing condition, which will be further studied in Section \ref{sc: d}. 

\begin{proposition}[\cite{Chernicketal1991}, Corollary 1.3]\label{T:Cor1.3} Let $\{X_n\}$ be a strictly stationary sequence of random variables such that for some $d\geq 1$ the conditions $D(u_n)$ and $D^{(d)}(u_n)$ hold for $u_n=u_n(\tau)$ for all $\tau>0$. Then the extremal index of $\{X_n\}$, $\theta$ exists  if and only if 
\begin{equation}\label{E:thetaCorr}
\lim_{n\to\infty}\P(M_{2,d}\leq u_n|X_1>u_n) = \theta,
\end{equation}
for all $\tau>0.$
\end{proposition}

Based on this limiting result, \cite{Hsing1993} has considered the empirical counterpart of the conditional probability in \eqref{E:thetaCorr} as  an estimator of $\theta$, assuming that $d$ is known. We shall follow this track, however first we resolve the problem of identifying a proper $d$. In practice, $d$ is unknown and moreover the $D^{(d)}(u_n)$ condition is not always satisfied, cf.\ examples in Section \ref{sec:simulation}. We shall make use of the stable tail dependence function to express $\theta$ and the $D^{(d)}(u_n)$ condition in the next subsection. Then, we will develop a statistical inference procedure to test if the $D^{(d)}(u_n)$ condition is fulfilled or not and to estimate a proper $d$. Last, an estimator of $\theta$ will be obtained using the estimated $d$.

\subsection[The link between $\theta$ and the STDF]{The link between $\theta$ and the stable tail dependence function}\label{sec: Delta_d}

Observe that  \eqref{E:thetaCorr} and the $D^{(d)}(u_n)$  condition are about the tail dependence structure of the random vectors $(X_1,\dots,X_d)$ and $(X_1,\dots,X_{r_n})$, respectively. 
To specify the notion of tail dependence, we introduce the stable tail dependence function $\ell_d$ for the random variables $(X_1,\ldots,X_d)$. If it exists, $\ell_d$ is defined as:
\begin{align}
\ell_d(x_1,\ldots,x_d) 
= \lim_{n\to\infty}n\P\left(X_1> u_n(x_1) \mbox{ or } \dots \mbox{ or } X_d> u_n(x_d) \right), 
 \label{E:ell}
\end{align}
for $(x_1,\dots,x_d)\in\R_+^d$ and $d\geq 1$. 
It is easy to see that the stable tail dependence function has homogeneity of order 1: for a positive constant $a$,
\begin{equation}
\ell_d(ax_1,\ldots,ax_d)=a\ell_d(x_1,\ldots,x_d). \label{eq: homo_l}
\end{equation}
See Chapter 6 in \cite{deHaanFerreira2006} for more properties of $\ell_d$. The tail process  introduced in \cite{BasrakSegers2009} is a commonly used tool for modeling the extremal dependence of a stationary sequence.  
Suppose that the tail process  of $(X_t)$ exists and we denote the tail process by $(Y_t)$. Then  $\ell_d$ also exists for any finite $d$ and moreover for $x_1\geq 1$,
$$
\ell_d(x_1,\ldots,x_d)=\P\left(Y_1>x_1^{\gamma}  \mbox{ or } \dots \mbox{ or } Y_d>x_d^\gamma\right),
$$
where $\gamma>0$ is the extreme value index of $X_1$. 

We denote by $\ell_d(\vc_d)$ the value of function $\ell_d$ evaluated at the $d$-dimensional unity vector $\vc_d:=(1,\dots,1)$, $d\geq 1$, and $\ell_0(1):=0$. The link between $\theta$ and $\ell_d$ is established in the following theorem.

\begin{theorem} \label{thm:d^*}
Assume that conditions $D(u_n)$ and $D^{(d)}(u_n)$ hold for some $d\geq 1$ and $\ell_{d}$ exists. Then, 
    \begin{equation}
	\theta = \ell_d(\vc_d)-\ell_{d-1}(\vc_{d-1})\in[0,1]. \label{E:defTheta}
    \end{equation}
Moreover, if $\ell_{s_0}$ exists for some $s_0\geq d$,  then there exists a positive integer $d^*\leq d$ such that 
for any $d^*\leq s\leq s_0$, $D^{(s)}(u_n)$ condition holds and 
$$
 \ell_s(\vc_s)-\ell_{s-1}(\vc_{s-1})=\theta;
$$
and, for any $1\leq s<d^*$,
$$
 \ell_s(\vc_s)-\ell_{s-1}(\vc_{s-1})>\theta.
$$
\end{theorem}
\paragraph*{Proof.}  By Proposition \ref{T:Cor1.3}, in order to prove \eqref{E:defTheta} it is sufficient to show that for any $u_n(\tau)$ such that 
$\lim_{n\rightarrow\infty}n\P(X_1>u_n(\tau))=\tau$, we have 
$\lim_{n\rightarrow \infty}\P(M_{2,d}\leq u_n|X_1>u_n)=\ell_d(\vc_d)-\ell_{d-1}(\vc_{d-1}).$
In fact,  this follows from,
\begin{align*}
 &\lim_{n\rightarrow \infty}\P(M_{2,d}\leq u_n(\tau)|X_1>u_n(\tau))\\
=&\lim_{n\rightarrow\infty}\frac{n}{\tau}\P(M_{2,d}\leq u_n(\tau), X_1>u_n(\tau))\\
=&\lim_{n\rightarrow\infty}\frac{n}{\tau}\left(\P(M_{2,d}\leq u_n(\tau))-\P(M_{1,d}\leq u_n(\tau))\right)\\
=&\lim_{n\rightarrow\infty}\frac{n}{\tau}\left(\P(M_{1,d}> u_n(\tau))-\P(M_{2,d}> u_n(\tau))\right)\\
=&\ell_d(\vc_d)-\ell_{d-1}(\vc_{d-1}).
\end{align*}
Under the assumption that $\ell_{s}$ exists, 
$$
\Delta(s):=\lim_{n\rightarrow\infty}\P(M_{2,s}\leq u_n|X_1>u_n)=\ell_s(\vc_s)-\ell_{s-1}(\vc_{s-1})
$$
is well defined for any finite $s$. Now since $\Delta(s)$ is a non-increasing function in $s$, we have, for any $s_0\geq d$, 
$$
\Delta(d)\geq \Delta(s_0)\geq \lim_{n\rightarrow\infty}\P(M_{2,r_n}\leq u_n|X_1>u_n).
$$
On the other hand, $D^{(d)}(u_n)$ condition is equivalent to the following equality:
\begin{equation}
\Delta(d)= \lim_{n\rightarrow\infty}\P(M_{2,r_n}\leq u_n|X_1>u_n).  \label{eq: Ddl}
\end{equation}
Therefore, 
 \begin{equation*}
\theta=\Delta(d)=\Delta(d+1)=\cdots=\Delta(s_0).
    \end{equation*}
Thus, $d^*$ exists and 
\begin{equation}
d^*=\min\{s: \Delta(s)=\theta\}.  \label{eq:d*}
\end{equation}
By the definition of $d^*$ and the monotonicity of $\Delta$, we have for any $1\leq s<d^*$, $\Delta(s)>\theta$.
\hfill  $\square$

\begin{remark}
If $d^*=1$, it implies that $D^{(1)}(u_n)$ condition holds and thus  $\theta=1$ by \eqref{E:defTheta} with $d=1$. 
And if $d^*\geq 2$, from the monotonicity of $\Delta(s)$ it follows that $\theta<1$. Moreover, \eqref{eq: Ddl} and \eqref{E:defTheta} motivate to 
check $D^{(d)}(u_n)$ condition and to obtain the value of $\theta$ by evaluating the function $\ell_s$. In Section \ref{sec:simulation}, we demonstrate the calculation for four models. 
\end{remark}

\subsection{Estimation of $\Delta(d)$}
In this subsection, we study an estimator of $\Delta(d)$ for some integer $d\in[d^*, s_0]$ and we shall establish the corresponding asymptotic normality. The estimator will be further used in testing $D^{(d)}(u_n)$ condition and in the estimation of $\theta$.

Recall the definition of $\Delta(d)=\lim_{n\rightarrow\infty}\P(M_{2,d}\leq u_n|X_1>u_n)$. 
First, we choose the intermediate quantile $F^{-1}(1-k/n)$ as the threshold $u_n$,  where 
 $k=k(n)$ is an intermediate sequence such that $k\rightarrow\infty$ and $n/k\rightarrow\infty$ as $n\rightarrow\infty$.  We shall estimate this threshold with 
$X_{n-k,n}$, where $\{X_{1,n}\leq X_{2,n}\leq \cdots \leq X_{n,n} \}$ are the order statistics of the sample. Then making use of the empirical probability measure, we have
\begin{align*}
\Delta(d)\approx &\frac{n}{k}\P(M_{2,d}\leq F^{-1}(1-k/n)< X_1)\\
\approx &\frac{n}{k}\cdot\frac{1}{n}\sum_{i=1}^{n-d+1}\mathds{1}\left\{M_{i+1,i+d-1}\leq X_{n-k,n}<  X_i  \right\}\\
=:&\widehat \Delta_n (d),
\end{align*}
where  $\mathds{1}$ denotes an indicator function. Note that for $d\geq d^*$, $\Delta(d)=\theta$. Thus, $\widehat \Delta_n (d)$ is also an estimator of $\theta$, which coincides with the estimator of $\theta$ in \cite{Hsing1993}. Next we derive the asymptotic normality of $\widehat \Delta_n(d)$ under some general mixing conditions.

Throughout we assume that $F$ is a continuous function. Let $U_i=1-F(X_i)$, $i=1,\ldots,n$, which become uniform random variables and let $\{U_{i,n}, i=1,\ldots,n\}$ be the order statistics. Observe that   $\mathds{1}\left\{M_{i+1,i+d-1}\leq X_{n-k,n}<  X_i  \right\}=
\mathds{1}\left\{  mU_{i+1,i+d-1}\geq U_{k,n} >U_i  \right\}$, where $mU_{i, j}:=\min_{i\leq t\leq j}U_j$ for $i\leq j$ and $mU_{i, j}:=1$ for $i>j$.
Therefore $\widehat \Delta_n (d)$ can be written as
\begin{equation}
\widehat \Delta_n (d)=\frac{1}{k}\sum_{i=1}^{n-d+1}\mathds{1}\left\{  mU_{i+1,i+d-1}\geq U_{k,n} >U_i  \right\}.  \label{eq: hat_theta}
\end{equation}
It is more convenient to formulate the conditions using 
$U_i$ than using $X_i$.

To obtain the asymptotic normality, we need a $\phi$-mixing conditions on the sequence  
$\{\mathds{1}\left\{U_i \leq \frac{k}{n}\right\}, i=1,\ldots, n\}$. Let $\mathcal{H}_l^s=\sigma(\mathds{1}\left\{U_i \leq \frac{k}{n}\right\}, l\leq i\leq s)$. Define
\begin{equation}
\phi_n(l)=\max_{s\geq 1}\sup_{A\in \mathcal{H}_i^s, B\in \mathcal{H}_{s+l}^n, \P(A)>0}|\P(B|A)-\P(B)|. \label{eq: phi}
\end{equation}
\begin{itemize}
\item[(A1)] There exist $r_n$ and $l_n$ such that $r_n\rightarrow\infty$, $\frac{r_n}{n}\rightarrow0$, $\frac{l_n}{r_n}\rightarrow 0$ and $\frac{n}{r_n}\phi_n(l_n)\rightarrow 0$.
\end{itemize}
Note that condition (A1) implies $D(u_n)$ condition and also the so called absolute regularity of the sequence, cf. \cite{Bradley2005}. We also need a strengthening condition of $D^d(u_n)$. 
\begin{itemize}
\item[(A2)] For $x, y\in [1/2, 3/2]$ and $r_n$ satisfying Condition (A1),
\begin{equation*}
\lim_{n\rightarrow \infty}\frac{n}{k}\sum_{i=d+1}^{r_n}\P(U_1<\frac{kx}{n}<mU_{2,d}, U_{i}<\frac{ky}{n})=0.  \label{eq:Dd2}
\end{equation*}
\end{itemize}
When $d=1$, Condition (A2) becomes the so called $D'(u_n)$ condition which implies $\theta=1$ and for $d=2$, it is equivalent to $D''(u_n)$ condition, 
cf. \cite{LeadbetterNandagopalan1989}. 

In order to obtain the asymptotic variance of $\widehat \Delta_n(d)$, we need the following two assumptions on the tail dependence structure of $(X_1,\ldots, X_{r_n})$.
\begin{itemize}
\item[(A3)] For $x, y\in [1/2, 3/2]$ and $r_n$ satisfying Condition (A1),
\begin{align*}
\lim_{n\rightarrow\infty}\frac{n}{k}\sum_{i=1}^{r_n}\P \left( U_1<\frac{kx}{n}, U_{i+1}<\frac{ky}{n}\right)
=\Lambda_1(x, y) \in [0,\infty).
\end{align*}
\item[(A4)] For $r_n$ satisfying Condition (A1),
\begin{align*}
\lim_{n\rightarrow\infty}\frac{n}{k}\sum_{i=1}^{r_n}\P \left( U_1<\frac{k}{n}, U_{i+1}<\frac{k}{n}<mU_{2+i,i+d}\right)
=\lambda_1 \in [0,\infty).
\end{align*}
\end{itemize}
The next condition makes sure that $\theta$ exists, cf.\ \eqref{E:defTheta}  and further it is used to control the bias of the limit distribution of the estimator. 
\begin{itemize}
\item[(A5)] There exists a $\rho>0$ such that for $j=d$ and $j=d-1$, as $t\rightarrow 0$,
\begin{equation}
 \frac{1}{t}\P\left(mU_{1,j}<t\right)-\ell_j(\vc_j)=O\left(t^\rho\right),
\end{equation}
where $\ell_j$ is defined in \eqref{E:ell}. 
\end{itemize}

\begin{theorem}  \label{thm:AN}
Assume that Condition $D(u_n)$ and Conditions (A1)-(A5) hold, $\sum_{i=r_n}^n\left(1-\frac{i}{n}\right)\phi_n(i)=o(1)$, $\frac{r_nk}{n}=o(1)$ and that $k=o\left(n^{2\rho/(2\rho+1)}\right)$. Then,
\begin{equation}
\sqrt{k}(\hat\Delta_n (d)-\theta)\overset{d}{\rightarrow} N(0,\sigma^2),
\end{equation}
where $\sigma^2=\theta(1-2\lambda_1)+\theta^2(2\Lambda_1(1,1)-1)$.
\end{theorem}

The proof of this theorem is provided in Section \ref{Sec:Proof}. Note that $\hat\Delta_n(d)$ is not yet an estimator of $\theta$ because in practice typically $d$ is unknown and  has to be estimated.  Our final estimator of $\theta$ is presented in Section \ref{sc: estimation_theta}.

\begin{remark}
\cite{Hsing1993} has already studied this estimator defined in \eqref{eq: hat_theta} and obtained its asymptotic normality under  $m$-dependence. 
\cite{WeissmanNovak1998} has established the asymptotic normality for runs and blocks estimators for a deterministic threshold, that is using $F^{-1}(1-k/n)$ (instead of $X_{n-k, n}$) as the threshold. Using a data-depending threshold and assuming a general mixing condition impose more technical challenges for proving the asymptotic result.  
\end{remark}

\subsection{Statistical inference on $D^{(d)}(u_n)$ condition} \label{sc: d}
We shall derive an inference procedure for testing $D^{(d)}(u_n)$ condition. The requirement of a suitable $d$ to satisfy the $D^{(d)}(u_n)$ condition is not limited 
to our estimation of $\theta$ alone. Some other methods of estimating $\theta$ also rely on the validity of the $D^{(d)}(u_n)$ condition for a specific value of $d$.
 In \cite{Suveges2007}, the author has proposed to check the $D^{(2)} (u_n)$ condition by looking at the empirical value of $\P(X_2\leq u_n< M_{3, r_n}|X_1>u_n)$. \cite{FerreiraFerreira2018} has used the same idea for checking $D^{(d)}(u_n)$ for a $d\geq2$. However, no formal consistency result has been established for such a procedure in both papers.

Throughout, we assume that $\theta$ exists and $\ell_{s_0}$ exists for some sufficiently large $s_0$.  For a given integer $d_0\geq 1$, we consider the following hypotheses:
$$
\begin{cases}
H_0:& \text{$D^{(d_0)}(u_n)$ condition holds.}\\
H_a:& \text{$D^{(d_0)}(u_n)$  condition does not  hold.}
\end{cases}
$$
From the proof of Theorem \ref{thm:d^*},  $\Delta(i)-\Delta(i+1)=0$ for $d_0\leq i <s_0$ under $H_0$. We investigate the asymptotic property of $\hat \Delta(i)-\hat \Delta(i+1)$ in the following theorem, based on which we will propose our test statistics. Define $\delta(s):=\Delta(s)-\Delta(s+1)$.

\begin{theorem}\label{thm:delta}
 Suppose that all the conditions of Theorem \ref{thm:AN} hold for $d=s+1$. Assume that Condition (A5) holds for $d=s$ and that the following limits exist and are finite.
 For any $x, y \in [1/2, 3/2]$,
\begin{align}
\sum_{i=s+1}^{r_n}\P\left(U_1<\frac{kx}{n}<mU_{2,s}, U_{i}<\frac{ky}{n}<mU_{i+1,i+s-1}\right) 
\rightarrow\Lambda_2(x, y) \in [0,\infty), \label{eq: condition_a2}
\end{align}
\begin{align*}
\sum_{i=1}^{r_n}\P \left( U_1<\frac{k}{n}, U_{i+1}<\frac{k}{n}<mU_{i+2,i+s}\right)
\rightarrow \tilde\lambda_1 \in [0,\infty),
\end{align*}
 \begin{align*}
\sum_{i=s+1}^{r_n}\P \left( U_1<\frac{k}{n}<mU_{2,s}, U_{i}<\frac{k}{n}<mU_{i+1,i+s}\right)
\rightarrow \lambda_2 \in [0,\infty),
\end{align*}
and
\begin{equation*}
\sum_{i=s+1}^{r_n}\P \left( U_1<\frac{k}{n}<mU_{2,s}, U_{i}<\frac{k}{n}\right)\rightarrow \lambda_3 \in [0,\infty).
\end{equation*}
Then, 
\begin{equation*}
\sqrt{k}\left( \widehat \Delta_n(s)-\widehat \Delta_n(s+1) -\delta(s) \right)\overset{d}{\rightarrow} N(0,\sigma^2_1),
\end{equation*}
where $\sigma^2_1=\delta^2(s)\left(2\Lambda_1(1,1)-1\right)-2\delta(s)\left(\tilde\lambda_1-\lambda_1+\lambda_3-\frac{1}{2}\right)+2\Lambda_2(1,1)-2\lambda_2.$
\end{theorem}
The proof of Theorem \ref{thm:delta} is provided in Section \ref{Sec:Proof}.

\begin{corollary}\label{cor:delta}
Assume that condition (A2) holds for $d=d^*$ and the conditions of Theorem \ref{thm:delta} hold for some $s\geq d^*$. Then 
\begin{equation}
\sqrt{k}\left( \widehat \Delta_n(s)-\widehat \Delta_n(s+1) \right)\overset{d}{\rightarrow} 0.  \label{eq:ss-1}
\end{equation}
 If $d^*\geq 2$ and the conditions of Theorem \ref{thm:delta} hold for $s=d^*-1$, then
\begin{equation}
\widehat \Delta_n(d^*-1)-\widehat \Delta_n(d^*) =\delta(d^*-1) +O_p\left(\frac{1}{\sqrt{k}}\right),   \label{eq:ss-2}
\end{equation}
where $\delta(d^*-1)>0$.
\end{corollary}
\paragraph*{Proof.} For $s\geq d^*$, $\delta(s)=\theta-\theta=0$. Moreover, Condition (A2) with $d=d^*$ implies that $\Lambda_2(1,1)=0$ and $\lambda_2=0$. Thus, $\sigma^2=0$  and \eqref{eq:ss-1} follows.  And, \eqref{eq:ss-2} is a direct consequence of Theorem \ref{thm:delta}  and the definition of $d^*$.

\hfill  $\square$\\
In view of \eqref{eq:ss-1}, we propose to accept $D^{(d_0)}(u_n)$ condition
if
\begin{equation}
\sqrt{k} \max_{d_0\leq s< s_0} \left(\widehat \Delta_n(s)-\hat\Delta_n(s+1)\right)<1,  \label{eq: Htest}
\end{equation}
where $s_0$ needs to be pre-specified. 
Because of the degenerate limit in \eqref{eq:ss-1}, there is no significance level in this testing procedure. So, it is not a standard hypothesis test. 
Under $H_a$, there are two scenarios. First, $d^*$ exists and $d^*>d_0$, so $D^{(d^*)}(u_n)$ holds. Then we have by \eqref{eq:ss-2},
\begin{align*}
\sqrt{k} \max_{d_0\leq s< s_0} \left(\widehat \Delta_n(s)-\hat\Delta_n(s+1)\right) 
\geq \sqrt{k} \left(\widehat \Delta_n(d^*-1)-\hat\Delta_n(d^*)\right) >>1.
\end{align*}
Under this situation, the power of this test tends to one asymptotically. 
The second scenario is that  $D^{(d)}(u_n)$ condition does not hold for any finite $d$. Then the conditions for deriving the asymptotic property of the test statistics in Corollary \ref{cor:delta} are not satisfied. We shall investigate the finite sample performance of this testing procedure for both scenarios in the simulation session.

\subsection{Estimation of $\theta$} \label{sc: estimation_theta}
Observe that $\hat\Delta_n(d)$ is a natural estimator of $\theta$ provided that $D^{(d)}(u_n)$ condition holds. The choice of $d$ is not unique: any $d^*\leq d\leq s_0$ would work. 
In practice, one can choose a $d$ using the hypothesis test proposed in Section \ref{sc: d}.  For the effectiveness of the estimation, we propose to use the smallest possible value, namely, $d^*$. It is difficult to prove that $d^*$ is the optimal choice from the theoretical point of view. From the simulation studies, we see that $d^*$ leads to the smallest mean squared error.

We need to estimate $d^*$ first.
Again following the idea that  $\Delta(s)-\Delta(s+1)=0$ for $d^*\leq s <s_0$, we propose to estimate $d^*$ by
\begin{equation}
\hat d^*(k)=\min \left\{s: \max_{s\leq i\leq s_0}\left(\widehat \Delta_n(i)-\hat\Delta_n(i+1)\right)<\frac{1}{\sqrt{k}}\right\}. \label{eq: hat_d}
\end{equation}
Provided that $s_0\geq d^*$ and the assumptions of Corollary \ref{cor:delta} hold, we have that (cf. Proposition \ref{prop:d*}) as $n\rightarrow\infty$,
$$
\P\left(\hat d^*(k)=d^*\right){\rightarrow} 1.
$$
Plugging $\hat d^*=\hat d^*(k)$ into \eqref{eq: hat_theta} yields our final estimator of $\theta$:
\begin{equation}
\hat \theta_n:=\widehat \Delta_n(\hat d^*)=\frac{1}{k}\sum_{i=1}^{n-\hat d^*+1}\mathds{1}\left\{  mU_{i+1,i+\hat d^*-1}\geq U_{k,n} >U_i  \right\}.  \label{eq: hat_thetas}
\end{equation}

Due to the consistency of $\hat d^*$, we conclude that the estimator $\hat \theta_n(\hat d^*)$ retains the asymptotic normality as shown in Theorem \ref{thm:AN} with $d=d^*$.

One need to specify before-hand two parameters, $k$ and $s_0$, for the test statistics and for the estimation of $\theta$. The parameter $s_0$ specifies the upper bound of the searching range and the specific choice is not crucial for the procedures. We recommend that a value  between 8 and 12 is a sensible choice in practice and is large enough because  typically one is interested in the  $D^{(d)}(u_n)$ condition for $d\leq 3$.  As for the parameter $k$, one can choose two different values of $k$ for $\hat d^*(k)$ and for  $\hat \theta_n$. In practice, we recommend to first estimate $d^*$ and then plug in the estimate to \eqref{eq: hat_thetas}. 
The optimal choice of $k$ depends on the convergence rate of the underlying model. In this paper, we don't investigate this from a theoretical perspective and rather show that our procedures work stably for a  range of $k$ in the simulation study. 

\section{Examples and simulation study} \label{sec:simulation}
To demonstrate the finite sample performance of our hypothesis test and our estimation of $\theta$, we consider five processes, among which four satisfy $D^{(d)}(u_n)$ condition for some $d$. For these four processes, we show how to compute $\theta$ and to validate  $D^{(d)}(u_n)$ condition via the stable tail dependence function $\ell_s$. 
In the sequel, we drop the subscript of the unit vector $\vc$ when it is clear from the context: $\ell_s(\vc)=\ell_s(\vc_s)$. In this section, $u_n$ is chosen such that $n\P(X>u_n)=1$.

\begin{itemize}
\item {\bf AR-N} This is a classical AR(1) model with Gaussian margin.  For $z\in (-1, 1)$, define
$$
X_i=zX_{i-1}+\epsilon_i, ~~ i\geq 1, 
$$
where $X_0\sim N(0, 1/(1-z^2))$ and $\epsilon_i$ i.i.d from $N(0,1)$.\\
Because of the asymptotic independence of multivariate normal random variables, we have $\ell_s(\vc)=s$, for $s\geq 1$. Thus \eqref{eq: Ddl} and therefore also $D^{(d)}(u_n)$ are satisfied for any $d\geq 1$. For this model, $\theta=1$ and $d^*=1$.\\

When $z=0$, this model becomes an independent sequence. We also include this special case for the simulation. 

\item {\bf Moving Maxima} Let $X_i=\max_{1\leq j\leq m}\epsilon_{j+i}$, for $i\geq 1$, where $m\geq 2$ is a fixed constant and $\epsilon_i$ an i.i.d sequence with $\P(\epsilon_i\leq x)=\exp\left(-\frac{1}{mx}\right)$. Then $(X_i)$ is a stationary sequence with marginal distribution $F(x)=\exp\left(-\frac{1}{x}\right)$ and for $s\geq 2$, 
    \begin{align*}
	\ell_s(\vc)=&\lim_{n\to\infty}n\P\left(\max_{1\leq i\leq s}X_i>u_n \right) \\
			=&\lim_{n\to\infty}n\left(1-F^{\frac{s+m-1}{m}}(u_n)\right)=\frac{s+m-1}{m}.
    \end{align*}
	Clearly, $D^{(d)}(u_n)$ are satisfied for any $d\geq 2$. For this model, $d^*=2$ and $\theta=\ell_2(\vc)-1=\frac{1}{m}$.


\item {\bf Max AR} For $z\in[0, 1)$, define $X_i=\max \{zX_{i-1}, \epsilon_i\}$, $i\geq 2$, where $X_1=\frac{\epsilon_1}{1-z}$ and 
$\epsilon_i$ i.i.d  with $\P(\epsilon_i\leq x)=\exp\left(-\frac{1}{x}\right)$.  Then for $s\geq 2$,
\begin{align*}
\ell_s(\vc)
=&\lim_{n\rightarrow \infty} n\left(1-\P\left( \max_{1\leq i\leq s}X_i\leq u_n\right)\right)\\
=&\lim_{n\rightarrow \infty} n\left(1-\P\left( \epsilon_1\leq (1-z)u_n, \max_{2\leq i\leq s}\epsilon_i\leq u_n\right)\right)\\
=&s-sz+z.
\end{align*}
Thus, $d^*=2$ and $\theta=1-z$.


\item {\bf AR-C} This is an AR(1) model with Cauchy margin. For $z\in (-1, 1)$, define
\begin{equation}
X_i=zX_{i-1}+\epsilon_i, i\geq 1, \label{eq:ARC}
\end{equation}
where $X_0$ has a standard Cauchy density $\frac{1}{\pi(1+x^2)}$ and $\epsilon_i$ i.i.d.\ with density $\frac{1-|z|}{\pi(x^2+(1-|z|)^2)}$.

By Proposition \ref{prop: ARC}, for $z\geq 0$, $d^*=2$ and $\theta=1-z$; and for $z<0$, $d^*=3$ and $\theta=1-|z|^2$.
\end{itemize}

Table \ref{tab: par} presents the parameters for each distribution that we have chosen for the simulation. 
\begin{table}[h]
\caption{Parameters of four models} \label{tab: par}
\centering
\begin{tabular}{l|ccccc}
 &IID &{AR-N} &Moving Maximum &MAX-AR  &{AR-C}\\
\hline
\hline \\[-.5ex]
$m$ or $z$ &-&  {$1/2$} &3 & {$1/2$}& {$-1/2$} \\[.2ex]
$d^*$&1  &1& 2& 2&3\\ [.2ex]
$\theta$ &1&1& $1/3$ &$1/2$ & $3/4$ 
\end{tabular}
\end{table}

Finally, we also consider an ARCH model which does not satisfy $D^{(d)}(u_n)$ condition for any finite $d$. The proof for this is given in Proposition \ref{prop: ARCH}.
\begin{itemize}

\item {\bf ARCH} Let
\begin{equation}
X_i=(2\times10^{-5}+0.7X^2_{i-1})^{1/2}\epsilon_i, ~~i\geq 1, \label{eq:ARCH}
\end{equation}
where $(\epsilon_i)$ are i.i.d. from $N(0, 1)$. For this model, $\theta=0.721$ (Table 3.2 in \cite{deHaanResnickRootzendeVries1989}).
\end{itemize}

\begin{table*}[h]
\centering
\caption{Acceptance (\%) of the hypothesis test for  $H_0:$ $D^{(d_0)}(u_n)$ condition holds. }\label{tab: Htest}
\begin{tabular}{l c c c c c c c}
\hline
 &\multicolumn{3}{c}{$k=50$} && \multicolumn{3}{c}{$k=100$} \\
Model & {$d_0=1$} & {$d_0=2$}& {$d_0=3$} && {$d_0=1$} & {$d_0=2$}& {$d_0=3$} \\
\hline
IID ($d^*=1$)&100 &100&100 &&100&100 &100\\
AR-N ($d^*=1$)&70.5 &100 &100 & &3.1 &99.9 &100\\
Moving Maxima ($d^*=2$)&0 &100 &100 & &0 &100 &100\\
MAX-AR ($d^*=2$)&0 &100 &100 & &0 &100 &100  \\
AR-C ($d^*=3$)&4.2 &3.3 &100 & &0 &0 &100  \\
ARCH &22.2 &97.6&100 & &0.4 &90.1 &100  \\
\hline
\end{tabular}
\end{table*}

In the simulations, we first look at the performance of the hypothesis test given in \eqref{eq: Htest} with $s_0=10$. We generate $n=5000$ observations from each distribution and compute the acceptance rate based on 1000 repetitions for $d_0= 1, 2, 3$. The result is reported in Table \ref{tab: Htest}.  For the IID and AR-N models, all the three null hypotheses are true because $d^*=1$. The test works ideally for the IID case while for the AR-N model it has a very high type II error for $d_0=1$, respectively, $29.5\%$ for $k=50$ and $96.9\%$ for $k=100$. However, for $d_0=2$ and 3, the type II error is (nearly) zero for both choices of $k$. For the Moving Maxima and MAX-AR models, the inference achieves the optimal accuracy: type I and type II errors are both zero. 
 For the AR model with the Cauchy margin, the procedure also performs optimally for $k=100$ and it has a type II error less than $5\%$ for $k=50$. 
For the ARCH model, which doesn't satisfy the $D^{(d)}(u_n)$ condition for any finite $d$, the procedure rejects $D^{(1)}(u_n)$ most of the time for $k=100$, however largely accepts $D^{(2)}(u_n)$ and $D^{(3)}(u_n)$ for both choices of $k$. This is not desired and it is probably due to the deviation from condition $D^{(d)}(u_n)$.

For the estimation of $\theta$, we set $s_0=10$ and use the same $k$ for both $\hat d^*$ and $\hat \theta_n$, which is convenient for implementation however not most efficient for the estimation.  In Section \ref{Sec:RealData}, we illustrate a way to first obtain $\hat d^*$ and then plot $\hat \theta_n$ against various $k$. 
To facilitate comparison, we also examine two existing methods  for estimating $\theta$.
\begin{itemize}
\item The pseudo maximum likelihood estimators based on the sliding blocks, $\hat \theta_n^{B,sl}$ defined on page 2314 in  \cite{BerghausBucher2018}. In this simulation, we take $n/k$ as the block length $r_n$ for this estimator.

\item The interval estimator from \cite{FerroSegers2003}, which is given by $\tilde \theta_n(u)$ on page 549 of that paper. We denote this estimator by $\hat \theta^{\rm{int}}$.

\end{itemize}

\begin{figure*}
\begin{center}
\includegraphics[width=0.85\textwidth]{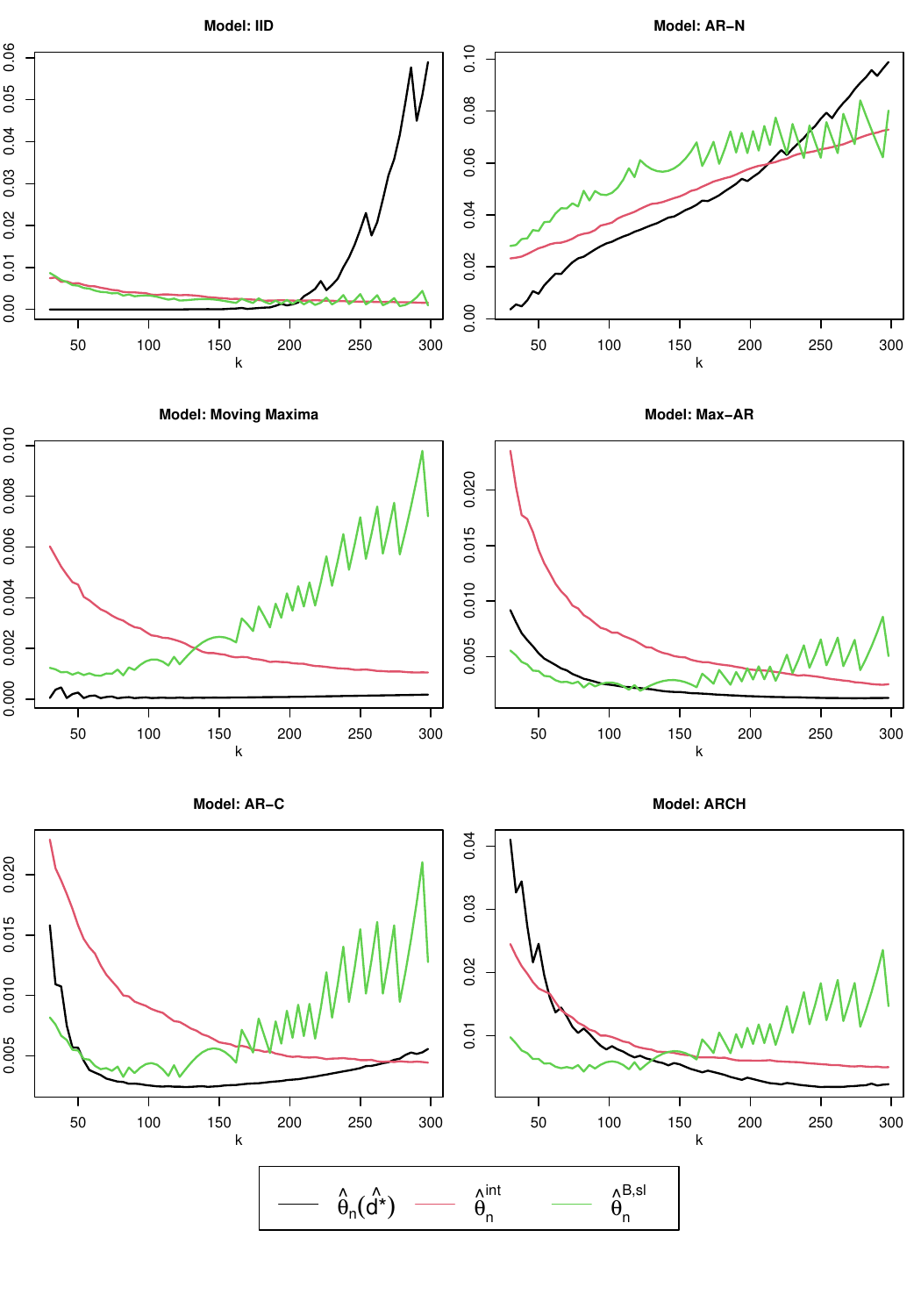}
\caption{The simulated MSE of three estimators for six distributions, where $\hat\theta(\hat d^*)$ is defined in \eqref{eq: hat_thetas}, $\hat \theta_n^{B,sl}$ defined in  \cite{BerghausBucher2018} and $\hat \theta^{\rm{int}}$ defined in \cite{FerroSegers2003}.}
\label{fig: MSE}
\end{center}
\end{figure*}

The simulated MSE is defined by \\$MSE(k)=\frac{1}{m}\sum_{i=1}^m (\hat \theta_i(k)-\theta)^2$, where $\hat\theta_i$ is the estimate based on the $i$-th repetition using one of the three methods and $m=1000$,  is plotted in Figure \ref{fig: MSE}. The dramatically increasing of MSE of our estimator for IID model is due to that $\hat d^*(k)$ often overestimates the true $d^*$ for a relatively large $k$. This can be largely improved by a two -step procedure: first estimate $d^*$ using a relative smaller $k$ and then take this estimate as an input in the estimate of $\theta$ without using the same $k$ for both. This procedure is illustrated in the next section.  For the other four  models that satisfy $D^{(d)}(u_n)$ condition,  our estimator outperforms the other two methods because it has the smallest MSE among the three for a sufficiently wide range of $k$. Even for the ARCH model, the minimum MSE of our estimator is also smaller than that of the other two methods. 
For the ARCH model, $\theta=\lim_{n\to\infty}\P(M_{2,r_n}\leq u_n|X_1>u_n)$. So $\theta$ can be well estimated by the runs estimator $\hat\Delta_n(r_n)$, with $\hat\Delta_n$ defined in \eqref{eq: hat_theta}. Our procedure leads to an estimator $\hat d^*$ such that the difference between $\widehat \Delta_n(\hat d^*)$ and $\widehat \Delta_n(s_0)$ is very small (cf. \eqref{eq: hat_d}). Therefore, for a finite sample our estimation of $d^*$ can be viewed as a selection procedure for the block length parameter of the runs estimator. 


\section{An application on summer daily temperature} \label{Sec:RealData}
We investigate the clustering of high temperature in summer by estimating the extremal index for data measured at three weather stations: de Bilt (N $52.11^\circ$, E $5.18 ^\circ$) in the Netherlands, Seville (N $37.4^\circ$, W $5.88 ^\circ$) in Spain, and Uccle (N $50.8^\circ$, E $4.35 ^\circ$) in Belgium. We consider the daily maximum temperatures in  July and August from 1951 to 1999 (\cite{rnoaa}). We choose this time interval because the data for Seville is not available prior to 1955 and the data for Uccle has many missing values after 1999. The sample size for each station is 3038. One gets more data if considering a single station. 

In terms of the measurement time, there is a natural gap between the data of two consecutive years, which are considered independent. To account for this feature, our estimator is adjusted as following: 
\begin{align}
 \widehat \Delta_n (d) 
=\frac{1}{k}\sum_{i=0}^{N-1}\sum_{j=1}^{L-d+1}\mathds{1}\left\{   U_{iL+j}<U_{k,n}\leq mU_{iL+j+1,iL+j+d-1} \right\}, \label{eq:tht_t}
\end{align}
where $N$ denotes the number of years, and $L$ denotes the number of observations for each year, 62 in this case.  

\begin{figure*}[t!]
\begin{center}
\includegraphics[width=0.95\textwidth]{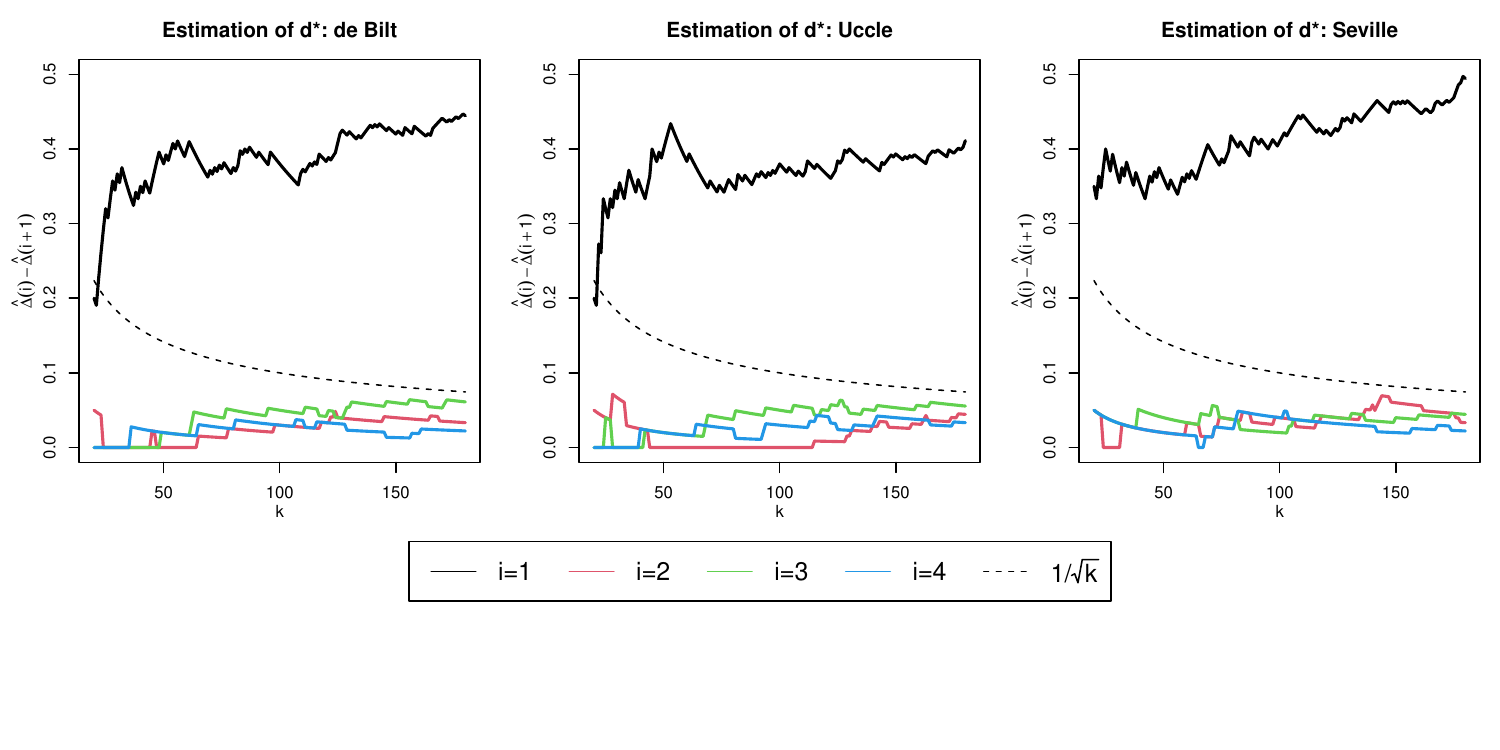}
\caption{The curves represent $\widehat \Delta_n(i)-\hat\Delta_n(i+1)$ (solid lines) and the threshold $1/\sqrt{k}$ (the dash lines), based on the summer temperature data observed in de Bilt, Uccle and Seville from 1951 to 1999. }
\label{fig: dtemp}
\end{center}
\end{figure*}

\begin{figure}
\begin{center}
\includegraphics[width=0.85\textwidth]{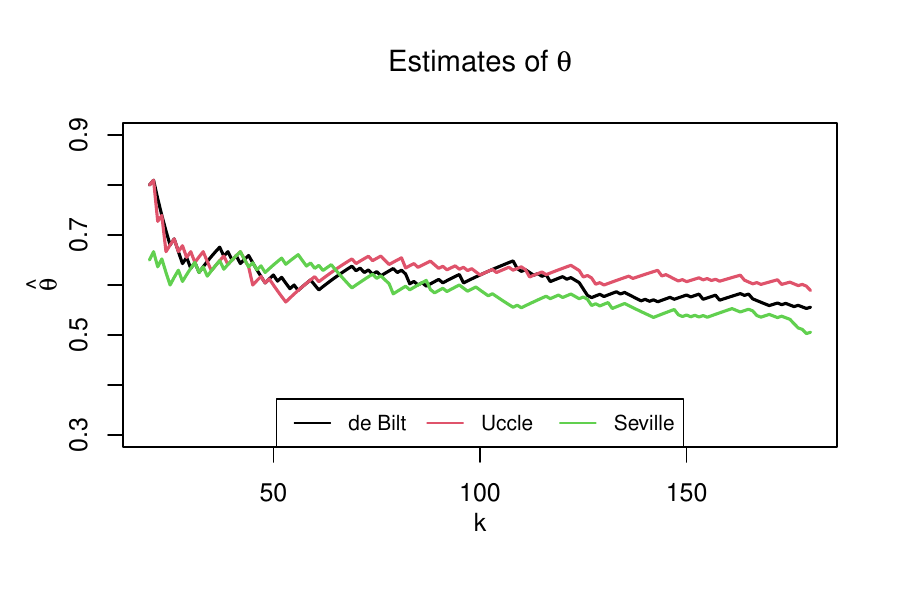}
\caption{Estimates of $\theta$ based on the summer temperature data observed in de Bilt, Uccle and Seville from 1951 to 1999. Note that $\hat \theta=\hat\Delta_n(2)$ given in \eqref{eq:tht_t}. }
\label{fig: dtht}
\end{center}
\end{figure}

We first estimate $d^*$ given by \eqref{eq: hat_d}. As illustrated in Figure \ref{fig: dtemp}, we have for all three stations $1-\hat\Delta_n(2)>1/\sqrt{k}$ and $\widehat \Delta_n(i)-\hat\Delta_n(i+1)< 1/\sqrt{k}$ for $i=2, 3,$ and 4. 
This provides a clear evidence that $\hat d^*=2$. The estimates of $\theta$ are plotted in Figure \ref{fig: dtht} by plugging $d=2$ into \eqref{eq:tht_t}. The estimates are stable for a wide range of $k$ and are remarkably close to each other for all three stations despite the different weather types.  In Table \ref{tab: tht_T}, we report the thresholds and the range of $\theta$ estimates for $k\in [50, 150]$.  Both the temperature thresholds and the estimates of $\theta$ in de Bilt resemble those in Uccle very closely. The thresholds in Seville are about \SI{10}{\celsius} higher than those in de Bilt and in Uccle. Yet,  the estimates of $\theta$ in Seville are close to the other two stations just with a slightly wider range.   With this case study, we conclude that the daily maximum temperature has a common extremal dependence structure in all three stations, that is, $D^{(2)}(u_n)$ condition is satisfied and $\theta$ is estimated to be around 0.6.

\begin{table*}[htbp]
\caption{ Thresholds and  the range of $\theta$ estimates for $k\in[50, 150]$} \label{tab: tht_T}
\centering
\begin{tabular}{l|lll}
 &De Bilt & Uccle & Seville\\
\hline
\hline \\[-.5ex]
$\left[X_{n-150,n}, ~~X_{n-50,n}\right]$& $\left[\SI{29.2}{\celsius} , \SI{31.3}{\celsius} \right] $& $\left[\SI{29.6}{\celsius} , \SI{32.0}{\celsius} \right]  $& $\left[\SI{40.8}{\celsius} , \SI{42.3}{\celsius} \right] $\\[.8ex]
$\left[\min_{k} \widehat \Delta_n (2), ~~ \max_{k} \widehat \Delta_n (2)\right]$& $[0.57, 0.64]$ & $[0.57, 0.66]$  & $[0.53, 0.67]$ 
\end{tabular}
\end{table*}

\section{Conclusion}\label{Sec:Conclusion}

Under the $D^{(d)}(u_n)$ condition, the extremal serial dependence of a stationary sequence is characterized by the quantity $\Delta(s)$ defined in Section \ref{sec: Delta_d}. Firstly,  $\Delta(s)$ can be used to formulate the $D^{(d)}(u_n)$ condition, capturing its essence.  Secondly, for an appropriate choice of $s$, typically not unique, $\theta$ can be equated to $\Delta(s)$. Thirdly, an elegant connection exists between $\Delta(s)$ and the stable tail dependence function, which enables the usage of the multivariate extreme value theory in the context of  the extremal serial dependence. 

The primary theoretical advancement presented in this paper is the establishment of the asymptotic normality of the empirical estimator $\widehat\Delta_n(s)$ for $\Delta(s)$.  Based on the asymptotic properties, we develop  hypothesis testing procedure for verifying the $D^{(d)}(u_n)$ condition. And it aids in estimating $d^*$, the smallest integer for which the $D^{(d)}(u_n)$ condition holds, which, in turn, is utilized in the estimation of $\theta$.


There are a few potential avenues for future research.. One natural yet challenging step would be developing asymptotic confidence intervals for $\theta$. However, a major obstacle to overcome is the dependence of the asymptotic variance of the estimator $\hat \theta_n$ (as described in Theorem \ref{thm:AN}) on unknown quantities such as $\lambda_1$ and $\Lambda(1,1)$. Thus, the estimation of these quantities becomes a prerequisite before constructing the confidence intervals.
Another promising direction for investigation is the utilization of estimators for the stable tail dependence function to estimate $\Delta(s)$ and, ultimately, $\theta$. However, a significant challenge lies in the fact that existing estimation methods for the stable tail dependence function are primarily designed for independent observations. Therefore, further research is necessary to adapt these estimators to account for serial dependence in the data.

\medskip

\noindent\textbf{Acknowledgment.}
 The author would like to acknowledge Andrea Krajina for the collaboration and fruitful
discussion in an early stage of this project, and Anja Jan{\ss}en for the discussion on the ARCH model. 



\bibliographystyle{imsart-number}
\bibliography{br}
\appendix

\section{Proof}\label{Sec:Proof}
To prove Theorems \ref{thm:AN} and \ref{thm:delta}, we first show two propositions. Define 
\begin{align}
\tilde\Delta_n(x,d)=\frac{1}{k}\sum_{i=1}^{n-d+1}\mathds{1}\left\{U_i<\frac{kx}{n}< mU_{i+1,i+d-1}\right\}.  \label{eq:def_ttht}
\end{align}
Note that $\tilde\Delta_n$ is a {\it pseudo} estimator because the $U_i$'s are not observable since $F$ is unknown. By the stationarity of the $U_i$'s,
\begin{align*}
\E\left(\tilde\Delta_n(x,d)\right) =\frac{n-d+1}{k}\P\left(U_1<\frac{kx}{n}< mU_{2,d}\right)
\rightarrow  \ell_d(x\vc_{d})-\ell_{d-1}(x\vc_{d-1}) =x\theta,
\end{align*}
by Condition (A5), the homogeneity of $\ell$ and \eqref{E:defTheta}.
And we also have,
\begin{align*}
\widehat \Delta_n (d)=\frac{1}{k}\sum_{i=1}^{n-d+1}\mathds{1}\left\{ U_i < U_{k,n}\leq mU_{i+1,i+d-1}\right\}
=\tilde\Delta_n\left(\frac{n}{k}U_{k, n} , d\right).
\end{align*}
Since $\frac{n}{k}U_{k, n}$ is expected to be around 1, we will first obtain the asymptotic properties for $\tilde\Delta_n(x,d)$ for $x\in[1/2, 3/2]$. 
Precisely, we shall prove the asymptotic normality of
$$
\sqrt{k}\left(\tilde \Delta_n(x,d)-\Delta_n(x,d)\right)=:\nu_n(x,d),~~~~ x\in[1/2, 3/2],
$$
where $\Delta_n(x,d)=\frac{n}{k}\P\left(U_1< \frac{kx}{n}< mU_{2,d} \right)$.



\begin{proposition} \label{prop:AN}
Under the conditions of Theorem \ref{thm:AN},
$$
\nu_n(x,d)\overset{d}{\rightarrow} W_d(x),
$$
in $D([1/2, 3/2])$, where $W_d$ is a continuous mean-zero Gaussian process with covariance structure
\begin{align*}
\E(W_d(x)W_d(y))
=\lim_{n\rightarrow\infty}\sum_{i=0}^{d-1}\frac{n}{k}\P\left(U_1<\frac{kx}{n}< mU_{2,d}, U_{1+i}<\frac{ky}{n}< mU_{2+i,d+i}\right),
\end{align*}
for $x\leq y$. In particular, $\E(W_d^2(x))=x\theta$. 
\end{proposition}
\noindent{\bf Proof of Proposition \ref{prop:AN}}\\
When there is no confusion, we drop the subscript $d$ from the notation of $W_d(x)$ and denote $\nu_n(x,d)=:\nu_n(x)$.
We prove convergence of the finite-dimensional distributions plus tightness. For the convergence of finite-dimensional distribution, it is sufficient to prove for each $m\in \N$ and for any  $x_i \in [1/2, 3/2]$ and $a_i\in \mathbf {R}$, $i=1,\ldots, m$,
$$
\sum_{i=1}^m  a_i \nu_n(x_i)\overset{d}{\rightarrow}\sum_{i=1}^m a_i W(x_i). 
$$
We show the proof for the case of $m=2$. For other cases, the proof is more tedious but can be done in the same way. Let $1/2\leq x\leq y \leq 3/2$, $I_i=\mathds{1}\left\{U_i<\frac{kx}{n}<mU_{i+1,i+d-1}\right\}$ and $J_i=\mathds{1}\left\{U_i<\frac{ky}{n}<mU_{i+1,i+d-1}\right\}$, $i=1,\ldots,n-d+1$. Then 
\begin{align*}
&a_1\nu_n(x)+a_2\nu_n(y)\\
=&\frac{1}{\sqrt{k}}\sum_{i=1}^{n-d+1}\left(a_1I_i+a_2J_i- \E(a_1I_i+a_2J_i)\right)+O\left(\frac{1}{\sqrt{k}}\right)\\
=:&\sum_{i=1}^{n-d+1}\xi_{i,n}+O\left(\frac{1}{\sqrt{k}}\right).
\end{align*}
 We apply the main theorem in \cite{Utev1990} to prove that $\sum_{i=1}^{n-d+1}\xi_{i,n}\overset{d}{\rightarrow}a_1W(x)+a_2W(y)$. We begin with computing the variance.
\begin{align*}
&\var\left(\sum_{i=1}^{n-d+1}\xi_{i,n}\right)\\
=&\frac{a^2_1}{k}\var\left(\sum_{i=1}^{n-d+1}I_i\right)
+\frac{a^2_2}{k}\var\left(\sum_{i=1}^{n-d+1}J_i\right)
+\frac{2a_1a_2}{k}\cov \left(\sum_{i=1}^{n-d+1}I_i, \sum_{i=1}^{n-d+1}J_i\right).
\end{align*}
The first two terms can be dealt with in the same way. Note that 
$$
\E(I_1)=\P\left(U_i<\frac{kx}{n}< mU_{i+1,i+d-1}\right)=O\left(\frac{k}{n}\right),
$$
and $\var (I_1)=\E(I_1)\left(1-O\left(\frac{k}{n}\right)\right)$. The same results hold for the $J_i$'s. Thus, by stationarity,
\begin{align*}
\var\left(\sum_{i=1}^{n-d+1}I_i\right)
=&(n-d+1)\var(I_1)+2\sum_{i<j}\cov(I_i, I_j)\\
=&(n-d+1)\E(I_1)\left(1-O\left(\frac{k}{n}\right)\right)
+2\sum_{i=1}^{n-d}(n-d+1-i)\cov(I_1, I_{1+i})\\
=&(n-d+1)\E(I_1)\left(1-O\left(\frac{k}{n}\right)\right)+o(k),
\end{align*}
by Lemmas \ref{lemma: cov} (i) and \ref{lemma: cov_mixing}. Thus,
\begin{equation*}
\frac{a_1^2}{k}\var\left(\sum_{i=1}^{n-d+1}I_i\right)
=a_1^2\frac{n}{k}\E(I_1)+o(1).
\end{equation*}
Similarly, one obtains that 
\begin{align*}
\frac{a^2_2}{k}\var\left(\sum_{i=1}^{n-d+1}J_i\right)=a_2^2\frac{n}{k}\E(J_1)+o(1).
\end{align*} 
As for the covariance term, we have that, again by stationarity and  Lemmas \ref{lemma: cov} (i) and \ref{lemma: cov_mixing},
\begin{align*}
&\cov \left(\sum_{i=1}^{n-d+1}I_i, \sum_{i=1}^{n-d+1}J_i\right)\\
=&\sum_{i=0}^{n-d+1} (n-d+1-i)\cov\left(I_1, J_{1+i}\right)
+\sum_{i=1}^{n-d+1} (n-d+1-i)\cov\left(I_{1+i}, J_{1}\right)\\
=&\sum_{i=0}^{d-1} (n-d+1-i)\cov\left(I_1, J_{1+i}\right)
+\sum_{i=1}^{d-1} (n-d+1-i)\cov\left(I_{1+i}, J_{1}\right)+o(k)\\
=&\sum_{i=0}^{d-1} (n-d+1-i)\E\left(I_1 J_{1+i}\right)
+\sum_{i=1}^{d-1} (n-d+1-i)E\left(I_{1+i}J_{1}\right)+O\left(\frac{k^2}{n}\right)+o(k).
\end{align*}
The last equality follows from that $\E(I_1)E(J_1)=O\left(\frac{n^2}{k^2}\right)$. Moreover, because of the disjointness causing by the fact that $x\leq y$, $\E(I_{i+1}J_1)=0$, for $i=1,\ldots,d-1$. Therefore, 
\begin{align*}
\frac{2a_1a_2}{k}\cov \left(\sum_{i=1}^{n-d+1}I_i, \sum_{i=1}^{n-d+1}J_i\right)
=2a_1a_2\frac{n}{k}\sum_{i=0}^{d-1} \E\left(I_1J_{1+i}\right)+o(1).
\end{align*}
We have shown that 
\begin{align*}
\lim_{n\rightarrow\infty} \var\left(\sum_{i=1}^{n-d+1}\xi_{i,n}\right)
=&\lim_{n\rightarrow\infty}\left(a_1^2\frac{n}{k}\E(I_1)+a_2^2\frac{n}{k}\E(J_1)+2a_1a_2\frac{n}{k}\sum_{i=0}^{d-1} \E\left(I_1, J_{1+i}\right)\right)\\
=&a_1^2\E(W^2(x))+a_2^2\E(W^2(y))+2a_1a_2\E(W(x)W(y)).
\end{align*}
Next, we check the Condition (2) in \cite{Utev1990}. Denote $\sigma^2_n=\var(\sum_{i=1}^{n-d+1}\xi_{i,n})$. Choosing 
$j_1=j_2=\cdots=1$, we have, for any $\epsilon>0$,
\begin{align}
\sigma_n^{-2}\sum_{i=1}^{n-d+1}\E(\xi_{i,n}^2\mathds{1}(|\xi_{i,n}|\geq \epsilon\sigma_n)) 
\leq &\sigma_n^{-3}n\epsilon^{-1}\E(|\xi_i|^3) \nonumber\\
=&\sigma_n^{-3}nk^{-3/2}\epsilon^{-1}\E\left(|a_1I_1+a_2J_1-\E(a_1I_1+a_2J_1)|^3\right) \nonumber\\
=&\sigma_n^{-3}nk^{-3/2}\epsilon^{-1}\E\left(|a_1I_1+a_2J_1-O(\frac{k}{n})|^3\right)\nonumber\\
=& c\sigma_n^{-3}nk^{-3/2}\epsilon^{-1} O\left(\frac{k}{n}\right)\rightarrow 0.  \label{eq: C2_Utev}
\end{align}
Thus, by the main theorem in \cite{Utev1990}, the central limit theorem holds for $\{\xi_{i,n}, i=1,\ldots, n-d+1\}$.

Next, we  prove that for any $\lambda>0$,
\begin{align}
\lim_{n\rightarrow \infty}\P\left(\sup_{\overset{|x-y|<\delta_n}{1/2\leq x<y\leq 2}}|\nu_n(x)-\nu_n(y)|\geq \lambda\right)=0,  \label{eq: tight_nu}
\end{align}
where for each $n$, $\delta_n$ is a constant and $\lim_{n\rightarrow\infty}\delta_n=0$.  Then by Theorem 1 in \cite{Aldous1978} and the explanation thereafter, $\nu_n$ is tight and each weak limit has a.s.\ continuous sample paths.  To prove \eqref{eq: tight_nu}, we decompose $\nu_n$ into two parts. Note that for any $x\in [1/2, 3/2]$,
$\mathds{1}\left\{U_i<\frac{kx}{n}< mU_{i+1,i+d-1}\right\}
=\mathds{1}\left\{mU_{i,i+d-1}<\frac{kx}{n}\right\}-\mathds{1}\left\{mU_{i+1,i+d-1}<\frac{kx}{n}\right\}.
$
 Thus, 
$$
\nu_n(x)=\nu_{1,n}(x)-\nu_{2,n}(x),
$$
where 
$\nu_{1,n}(x)=\frac{1}{\sqrt{k}}\sum_{i=1}^{n-d+1} (\mathds{1}\{mU_{i,i+d-1}<\frac{kx}{n}\}-\P(mU_{i,i+d-1}<\frac{kx}{n}))$,
and $\nu_{2,n}(x)=\frac{1}{\sqrt{k}}\sum_{i=1}^{n-d+1} (\mathds{1}\{mU_{i+1,i+d-1}<\frac{kx}{n}\}-P(mU_{i+1,i+d-1}<\frac{kx}{n}))$.

To prove \eqref{eq: tight_nu}, it is sufficient to prove  the tightness condition  for $\nu_{1,n}(x)$ and $\nu_{2,n}(x)$, respectively. We demonstrate the proof for $\nu_{1,n}(x)$. 
Let $t_n=\lfloor\frac{n-d+1}{2r_n}\rfloor$. We split the sum into $2t_n$ blocks of length $r_n$ and a remaining block of length less than $2r_n$. 
To simplify the notation, we denote $M_i=mU_{i,i+d-1}$. Precisely,

\begin{align*}
\nu_{1,n}(x)
=&\frac{1}{\sqrt{k}}\sum_{i=1}^{n-d+1}\left( \mathds{1}\left\{M_{i}<\frac{kx}{n}\right\}-\P\left(M_{i}<\frac{kx}{n}\right)\right)\\
=&\frac{1}{\sqrt{k}}\sum_{i=0}^{t_n-1}\sum_{j=1}^{r_n}\left(     \mathds{1}\left\{M_{2ir_n+j}<\frac{kx}{n}\right\} 
-\P\left(M_{2ir_n+j}<\frac{kx}{n}\right)   \right.\\
&+\left.     \mathds{1}\left\{M_{(2i+1)r_n+j}<\frac{kx}{n}\right\}-\P\left(M_{(2i+1)r_n+j}<\frac{kx}{n}\right)      \right)\\
&+\frac{1}{\sqrt{k}}\sum_{i=2t_nr_n+1}^{n-d+1}\left(     \mathds{1}\left\{M_{i}<\frac{kx}{n}\right\}-\P\left(M_{i}<\frac{kx}{n}\right)      \right)\\
=:&\mu_{1,n}(x)+\mu_{2,n}(x)+\mu_{3,n}(x).
\end{align*}

Define $\tilde\mu_n(x)=\frac{1}{\sqrt{k}}\sum_{i=1}^{t_nr_n}(\mathds{1}\{\tilde M_i<\frac{kx}{n}\}-\P(\tilde M_i<\frac{kx}{n}))$, where for  any $i=1,\ldots, t_n$, 
$$\{\tilde M_{r_n(i-1)+i},\ldots, \tilde M_{ir_n}\}\overset{d}{=}\{ M_{1},\ldots, M_{r_n}\},$$ 
and  $\{\tilde M_{r_n(i-1)+j}, j=1,\ldots, r_n\}_{i=1}^{t_n}$ are $t_n$ independent blocks. So the $\tilde M_i$'s form a special $r_n$-dependent array for each $n$, which is not a strictly stationary sequence. We first apply a fluctuation inequality for $m$-dependent arrays given by Theorem 4.1 in \cite{EinmahlRuymgaart2000}  to prove the tightness of $\tilde\mu_n(x)$. Then, the tightness of $\mu_{1,n}$ and $\mu_{2,n}$ follows from the bounded variation distance between $\tilde\mu_n$ and $\mu_{1,n}$ and between $\tilde \mu_n$ and $\mu_{2,n}$, respectively.

For each $n$, let $q={r_n^{1+\epsilon}}$, where $\epsilon$ is some positive number such that $q/\sqrt{k_n}\rightarrow 0$. Define $I_i=[\frac{1}{2}+\frac{i}{q}, \frac{1}{2}+\frac{i+1}{q}],  ~~i=0,\ldots, q-1$. Choose 
$\delta_n=\frac{1}{q}$. For any $x, y \in [1/2, 3/2]$ and $|x-y|<\delta_n$, there exists an $i \in \{1, \ldots, q-1\}$ such that 
$$
\left|x-\frac{1}{2}-\frac{i}{q}\right|<\frac{1}{q} ~~~~~\text{ and }~~~~\left|y-\frac{1}{2}-\frac{i}{q}\right|<\frac{1}{q}. 
$$
Thus, for any $\lambda>0$,

\begin{align*}
&\P\left(\sup_{\overset{|x-y|<\delta_n}{1/2\leq x<y\leq 2}}|\tilde \mu_n(x)-\tilde\mu_n(y)|\geq \lambda\right)\\
\leq &\P\left(\max_{1\leq i\leq q-1} \sup_{\overset{|x-1/2-{i}/{q}|<{1}/{q}}{|y-1/2-{i}/{q}|<{1}/{q}}} ( |  \tilde \mu_n(x)-\tilde\mu_n(\frac{1}{2}+\frac{i}{q})| 
+ |\tilde \mu_n(\frac{1}{2}+\frac{i}{q})-\tilde\mu_n(y)|  ) \geq \lambda\right)\\
\leq &2\P\left( \max_{1\leq i\leq q-1} \sup_{x\in I_i}\left|\tilde \mu_n(x)-\tilde\mu_n(\frac{1}{2}+\frac{i}{q})\right| \geq \frac{\lambda}{2}  \right)
+2\P\left( \max_{1\leq i\leq q-1} \sup_{x\in I_i}\left|\tilde \mu_n(x)-\tilde\mu_n(\frac{1}{2}+\frac{i+1}{q})\right| \geq \frac{\lambda}{2}  \right)\\
\leq & 2\sum_{i=1}^{q-1}P\left( \sup_{x\in I_i}\left|\tilde \mu_n(x)-\tilde\mu_n(\frac{1}{2}+\frac{i}{q})\right| \geq \frac{\lambda}{2}  \right)
+2\sum_{i=1}^{q-1}P\left( \sup_{x\in I_i}\left|\tilde \mu_n(x)-\tilde\mu_n(\frac{1}{2}+\frac{i+1}{q})\right| \geq \frac{\lambda}{2}  \right)\\
\leq  & 4\sum_{i=1}^{q-1}P\left( \sup_{x, y\in I_i}\left|\tilde \mu_n(x)-\tilde\mu_n(y)\right| \geq \frac{\lambda}{2}  \right)
=:4\sum_{i=1}^{q-1}P_i.
\end{align*}

Next, we apply (4.4) in \cite{EinmahlRuymgaart2000} to bound $P_i$. Define
$$
\Gamma_n (x)=\frac{1}{t_nr_n}\sum_{i=1}^{t_nr_n}\mathds{1}\left\{\tilde M_i<x\right\}-\P(\tilde M_i<x),
$$
which plays the role of $\Delta_n$ in \cite{EinmahlRuymgaart2000}.
Then, by taking, in the notation of that paper, $\epsilon =\frac{1}{2}$ and $m=r_n$, we have
\begin{align*}
P_i=&\P\left(\sup_{\frac{1}{2}+\frac{i}{q}\leq x<y\leq \frac{1}{2}+\frac{i+1}{q}}|\tilde \mu_n(x)-\tilde\mu_n(y)|\geq \lambda\right)\\
=&\P\left(\frac{r_nt_n}{\sqrt{k}}\sup_{\frac{1}{2}+\frac{i}{q}\leq x<y\leq \frac{1}{2}+\frac{i+1}{q}}|\Gamma_n\left(\frac{kx}{n}\right)-\Gamma_n\left(\frac{ky}{n}\right)|\geq \lambda\right)\\
=&\P\left(\sup_{\frac{k}{n}(\frac{1}{2}+\frac{i}{q})\leq a<b\leq\frac{k}{n}(\frac{1}{2}+\frac{i+1}{q})}|\Gamma_n\left(a\right)-\Gamma_n\left(b\right)|\geq \frac{\sqrt{k}\lambda}{r_nt_n}\right)\\
\leq& c_1\exp\left(\frac{-r_nt_n\frac{k\lambda^2}{r_n^2t_n^2}}{4r_np_i}\psi\left(\frac{\sqrt{r_nt_n}  \frac{\sqrt{k}\lambda}{r_nt_n}}{\sqrt{r_nt_n}p_i} \right)\right)
= c_1\exp\left(-\frac{k\lambda^2}{4r_n^2t_np_i}\psi\left(\frac{\sqrt{k}\lambda}{r_nt_np_i}   \right)\right)
\end{align*}
where $p_i=\P\left(\frac{k}{n}(\frac{1}{2}+\frac{i}{q})<mU_{1,d}\leq \frac{k}{n}(\frac{1}{2}+\frac{i+1}{q})\right)$ and $\psi$ is a continuous and decreasing function such that $\psi(0)=1$. Observe that 
\begin{align*}
\sup_{0\leq i\leq q-1}\left| \frac{n}{k}p_i-\left(\ell_d\left(\left(\frac{1}{2}+\frac{i+1}{q}\right)\vc_d\right)- 
 \ell_d\left(\left(\frac{1}{2}+\frac{i}{q}\right)\vc_d\right)\right)\right|=o\left(\left(\frac{n}{k}\right)^\tau\right),
\end{align*}
and $\ell_d\left(\left(\frac{1}{2}+\frac{i+1}{q}\right)\vc_d\right)-\ell_d\left(\left(\frac{1}{2}+\frac{i}{q}\right)\vc_d\right)=\frac{1}{q}\ell_d(\vc_d)$ by the homogeneity of $\ell_d$. Thus, for $n$ large enough, $\frac{k}{2nq}\leq p_i\leq \frac{2dk}{nq}$, uniformly in $i$, due to the fact that $\ell_d(\vc_d)\in [1, d]$.

Then by the choice of $q$ and that $n-r_n\leq 2 r_nt_n\leq n$,
$$\frac{k}{4r_n^2t_np_i}\geq \frac{k}{2r_nnp_i}\geq \frac{m}{4dr_n}=\frac{1}{4d}r_n^\epsilon,
$$
and $$ \psi\left(\frac{\sqrt{k}\lambda}{r_nt_np_i}   \right)\geq \psi\left(\frac{3\sqrt{k}\lambda}{n\frac{k}{2nq}} \right)
=\psi\left(\frac{6\lambda q}{\sqrt{k}}  \right)\rightarrow \psi(0)=1.  $$
 Thus,
\begin{align*}
\P\left(\sup_{\overset{|x-y|<\delta_n}{1/2\leq x<y\leq 2}}|\tilde \mu_n(x)-\tilde\mu_n(y)|\geq \lambda\right)
\leq 4q\exp(-c_2r_n^\epsilon)=4r_n^{1+\epsilon}\exp(-c_2r_n^\epsilon)\rightarrow 0.
\end{align*}
So the tightness of $\tilde \mu_n$ follows from the  tightness criterion by Theorem 1 in \cite{Aldous1978}.

  By Lemma 2 in \cite{Eberlein1984}, 
\begin{align*}
&\parallel \Omega\left(\{\tilde M_{r_n(i-1)+1},\ldots, \tilde M_{ir_n}\}_{i=1}^{t_n} \right) 
 -\Omega\left(\{ M_{r_n(i-1)+1},\ldots, M_{ir_n}\}_{i=1}^{t_n} \right)    \parallel \\
=&\parallel \bigotimes_{i=1} ^{t_n}\Omega\left( M_{r_n(i-1)+1},\ldots, M_{ir_n} \right) 
-\Omega\left(\{ M_{r_n(i-1)+1},\ldots, M_{ir_n}\}_{i=1}^{t_n} \right) \parallel \\
\leq &\phi_n(r_n)t_n\leq \beta(r_n) \frac{n}{r_n}\rightarrow 0,
\end{align*}
by the absolutely regular assumption on the sequence, and the condition that $\beta(r_n) \frac{n}{r_n}\rightarrow 0$, where $\Omega (X)$ denotes the distribution of $X$.
Thus,  we obtain that for $i=1, 2$,
\begin{align}
\P\left(\sup_{\overset{|x-y|<\delta_n}{1/2\leq x<y\leq 2}}| \mu_{i,n}(x)-\mu_{i,n}(y)|\geq \lambda\right)\rightarrow 0.
\end{align}

To prove the tightness of $\nu_{1,n}$, it is remaining to show that $\sup_{1/2\leq x\leq 3/2}|\mu_{3, n} (x)|\overset{P}{\rightarrow} 0$. Note that by the definition of $t_n$, the number of summands in $\mu_{3, n}$ is bounded by $2r_n$. 
\begin{align*}
\E\left(\sup_{1/2\leq x\leq 3/2}|\mu_{3, n} (x) \right)
\leq &\E\left( \frac{1}{\sqrt{k}}\sum_{i=2t_nr_n+1}^{n-d+1}\left(     \mathds{1}\left\{M_{i}<\frac{3k}{2n}\right\}+\P\left(M_{i}<\frac{3k}{2n}\right)      \right) \right)\\
\leq & \frac{4r_n}{\sqrt{k}}\cdot\frac{3k}{2n}\rightarrow 0,
\end{align*} 
by the assumption that $\frac{r_n\sqrt{k}}{n}\rightarrow 0$.

 \hfill  $\square$
 
By a similar but simpler proof, we can obtain the convergence for the empirical tail process: $\nu_n(x, 1)=\sqrt{k}\left(\frac{1}{k}\sum_{i=1}^n \mathds{1}\left(U_i<\frac{kx}{n}\right)-x\right)$. 
\begin{proposition} \label{prop: marginal_AN}
Under the conditions of Theorem \ref{thm:AN},
\begin{equation*}
\nu_n(x, 1)\overset{d}{\rightarrow} \tilde W(x),
\end{equation*}
in $D([1/2, 3/2])$, where $\tilde W$ is a continuous mean-zero Gaussian process with covariance structure\\
$\E(\tilde W(x) \tilde W(y))=\min (x, y)+\Lambda_1(x, y)+\Lambda_1(y, x)$. 
\end{proposition}
Note that this result is generally different from Proposition \ref{prop:AN} with $d=1$. The difference lies in the covariance structure. The two coincide with each other only when  Condition (A2)   holds with $d=1$, which implies $D^{(1)}(u_n)$ condition holds. 

 Recall that $\hat \theta_n=\tilde\Delta_n\left(\frac{n}{k}U_{k, n},d\right)$. We now derive asymptotic property for $\frac{n}{k}U_{k, n}=:e_n$.
First, Proposition \ref{prop: marginal_AN} implies that, by Theorem A.0.1 and Lemma A.0.2 in \cite{deHaanFerreira2006},
\begin{equation}
\sqrt{k}\left(e_n-1\right)\overset{d}{\rightarrow} -\tilde W(1).  \label{eq: quantileprocess}
\end{equation}
In particular, one has $e_n\overset{p}{\rightarrow} 1$. Further, by the continuity of $\tilde W$, 
\begin{align}
&\left\lvert  \sqrt{k}\left(e_n-1\right)+ \nu_n(1, 1) \right\rvert  \nonumber\\
=&\left\lvert  \nu_n(e_n, 1) - \nu_n(1, 1) \right\rvert     \nonumber\\
\leq& \left\lvert  \nu_n(e_n, 1) - \tilde W(e_n)\right\rvert +   \left\lvert  \nu_n(1, 1) - \tilde W(1)\right\rvert  
+ \left\lvert  \tilde W(1) - \tilde W(e_n)\right\rvert   
=o_p(1).    \label{eq: quantileprocess2}
\end{align}

\noindent{\bf Proof for Theorem \ref{thm:AN} }\\
For a convenient presentation, all the processes involved in the proof are defined on the
same probability space, via the Skorohod construction. We use the same notation,
though they are only equal in distribution to the original ones. 
We start with the following decomposition: by the definitions of $\tilde\Delta_n$ and $\nu_n(x, d)$
\begin{align*}
\sqrt{k} \left(\widehat \Delta_n(d)-\theta\right)
=\sqrt{k} \left(\tilde\Delta_n\left(e_n,d\right)-\theta\right)
=\nu_n(e_n, d)
+\sqrt{k}\left( \Delta_n\left(e_n,d\right) -\theta\right). 
\end{align*}
We denote $\sqrt{k}\left( \Delta_n\left(e_n,d\right) -\theta\right)=:I_n$. It is sufficient to prove the following three limit relations, from which the theorem follows. 
\begin{align}
\left| \nu_n(e_n, d)-\nu_n(1, d)\right|  \overset{p}{\rightarrow}0, \label{eq: i1}
\end{align}
\begin{align}
\left|I_n+\theta \nu_n(1,1)\right|  \overset{p}{\rightarrow}0, \label{eq: i2}
\end{align}
and 
\begin{align}
 \nu_n(1, d)-\theta \nu_n(1,1) \overset{d}{\rightarrow} N(0, \sigma^2), \label{eq: i12}
\end{align}
where $\sigma^2$ is defined in Theorem \ref{thm:AN}.

For \eqref{eq: i1}, we have 
\begin{align*}
&\left|  \nu_n(e_n, d)-\nu_n(1, d) \right|\\
=&\left|\left(\nu_n(e_n, d)- W_d\left(e_n\right)\right)+\left(W_d(1)-  \nu_n(1, d)\right) 
 +\left(W_d\left(e_n\right)-W_d(1)\right) \right|\\
\leq &2\sup_{1/2\leq x\leq 2}\left| \nu_n(x,d)-W_d(x)\right| +\left|W_d\left(e_n\right)-W_d(1)   \right|  \overset{p}{\rightarrow}0,
\end{align*}
by Proposition \ref{prop:AN} and the continuous mapping theorem,   continuity of $W_d$  and $e_n \overset{p}{\rightarrow} 1$. 

Next, we deal with $I_n$. Define $\ell_{d,n}(x\vc_d)=\frac{n}{k}\P\left( mU_{1,d}\leq \frac{kx}{n} \right)$. Then $\lim_{n\rightarrow\infty}\ell_{d,n}(x\vc_d)=\ell_d(x\vc_d)$ and $\Delta_n(x,d)=\ell_{d,n}(x\vc_d)-\ell_{d-1,n}(x\vc_{d-1})$.  Therefore, 
\begin{align*}
I_n=&\sqrt{k}\left( \ell_{d,n}\left(e_n\vc_d\right)-\ell_{d-1,n}\left(e_n\vc_{d-1}\right)- \theta\right)\\
=&\sqrt{k}\left(   \ell_{d,n}\left(e_n\vc_d\right)- \ell_{d}\left(e_n\vc_d\right) 
 -\left( \ell_{d-1,n}\left(e_n\vc_{d-1}\right)- \ell_{d-1}\left(e_n\vc_{d-1}\right)\right)\right)\\
&+\sqrt{k}\left( \ell_{d}\left(e_n\vc_d\right)- \left(  \ell_{d-1}\left(e_n\vc_{d-1}\right)- \theta\right)\right)\\
=:&I_{n1}+I_{n2}.
\end{align*}
By Assumption (A5), $k=o\left(n^{2\rho/(2\rho+1)}\right)$ and  $e_n \overset{p}{\rightarrow} 1$, $I_{n1}\overset{p}{\rightarrow}0$.
 By the homogeneity of $\ell_d$, \eqref{E:defTheta} and \eqref{eq: quantileprocess2},
\begin{align*}
I_{n2}=\theta\sqrt{k}\left(e_n-1\right)=-\theta \nu_n(1,1)+o_p(1).
\end{align*}
Thus, \eqref{eq: i2} is proved.


Last, to prove \eqref{eq: i12}, we denote $T_i=\mathds{1}\{U_i<\frac{k}{n}< mU_{i+1,i+d-1}\}-\theta \mathds{1}\{U_i<\frac{k}{n}\}$, then   $ \nu_n(1, d)-\theta \nu_n(1,1)=\frac{1}{\sqrt{k}}\sum_{i=1}^{n-d+1}\left(T_i- \E(T_1)\right)+O_p(d/\sqrt{k})$.
 We shall apply the main Theorem in \cite{Utev1990} to prove the central limit theorem for $\sum_{i=1}^{n-d+1}\xi_{i,n}$, where $\xi_{i,n}=\frac{1}{\sqrt{k}}\left(T_i- \E(T_1)\right)$. We begin with the variance: 
\begin{align*}
\var\left(\sum_{i=1}^{n-d+1}\xi_{i,n} \right)
=\frac{n-d+1}{k}\var(T_i)+\frac{2}{k}\sum_{i=1}^{n-d+1}(n-i)\cov(T_1, T_{1+i})
\end{align*}
First, $\frac{n}{k}\E(T_1)\rightarrow \theta-\theta=0$. Thus, 
\begin{align*}
&\frac{n}{k}\var(T_i)=\frac{n}{k}\E(T^2_1)+o(1)\\
=& \frac{n}{k}\left(\P \left(U_1<\frac{k}{n}<mU_{2,d}\right) +\theta^2 \P\left(U_1<\frac{k}{n}\right) 
-2\theta\P \left(U_1<\frac{k}{n}<mU_{2,d}\right) \right)+o(1)\\
\rightarrow &\theta-\theta^2.
\end{align*}
Second, by the mixing condition,\\ $\frac{2}{k}\sum_{i=r_n+1}^{n-d+1}(n-i)\cov(T_1, T_{1+i})=o(1)$. 
And  by Conditions (A2), (A3) and (A4),
\begin{align*}
&\frac{2}{k}\sum_{i=1}^{r_n}(n-i)\cov(T_1, T_{1+i})
=\frac{2n}{k}\sum_{i=1}^{r_n}\left(1-\frac{i}{n}\right)\E(T_1T_{i+1})+o(r_nk/n)\\
=&\frac{2n}{k}\sum_{i=1}^{r_n}\theta^2\P\left(U_1<\frac{k}{n}, U_{1+i}<\frac{k}{n}\right)
-\frac{2n}{k}\sum_{i=1}^{r_n}\theta \P\left(U_1<\frac{k}{n}, U_{1+i}<\frac{k}{n}<mU_{i+2,i+d}\right)\\
&-\frac{2n}{k}\sum_{i=1}^{r_n}\theta \P\left(U_1<\frac{k}{n}< mU_{2,d}, U_{i+1}<\frac{k}{n}\right)\\
&+\frac{2n}{k}\sum_{i=1}^{r_n}\P\left(U_1<\frac{k}{n}< mU_{2,d}, U_{i+1}<\frac{k}{n}<mU_{i+2,i+d}\right)
+o(r_nk/n)\\
\rightarrow &2\theta^2\Lambda_1(1,1)-2\theta\lambda_1.
\end{align*}
And the Condition (2) in \cite{Utev1990} follows from the same argument as that for \eqref{eq: C2_Utev}. Therefore,
\eqref{eq: i12} is proved.

 \hfill  $\square$

\noindent{\bf Proof for Theorem \ref{thm:delta}}\\  
Because the conditions of Theorem \ref{thm:AN} hold for $d=s+1$, $\Delta(s+1)=\theta$ and  the result in Proposition \ref{prop:AN} holds for $d=s+1$. 

We also need a similar convergence result for $\sqrt{k}(\tilde \Delta_n(x,s)-\Delta_n(x,s))$. Clearly, the covariance of the limit is not necessary in the same form  because $D^{(s)}$ condition is not guaranteed. However the conditions in Theorem \ref{thm:delta} makes sure that the covariance exists and the tightness of the process still holds for $d=s$. 
Follow the same line as in the proof for Proposition \ref{prop:AN}, we have
$$
\sqrt{k}\left(\tilde \Delta_n(x,s)-\Delta_n(x,s)\right)\overset{d}{\rightarrow} V(x),
$$
in $D([1/2, 3/2])$, where $V$ is a continuous mean-zero Gaussian process with covariance structure
\begin{align*}
\E(V(x)V(y))
=&\lim_{n\rightarrow\infty}\frac{n}{k}\sum_{i=0}^{s-1}\P\left(U_1<\frac{kx}{n}< mU_{2,s}, U_{1+i}<\frac{ky}{n}< mU_{2+i,s+i}\right)\\
&+\Lambda_2(x,y)+\Lambda_2(y,x),
\end{align*}
for $x\leq y$. Note that if Condition (A2) holds for $d=s$, then $\Lambda_2(x,y)=0$ and $V=^d W_{s}$. 
It is clear that
\begin{equation}
\E(V^2(1))=\Delta(s)+2\Lambda_2(1,1).  \label{eq: var_v1}
\end{equation}

Then the theorem can be proved in the same manner as that for Theorem \ref{thm:AN}. Recall that $\delta(x)=\Delta(s)-\Delta(s+1)$. By the same argument for 
 \eqref{eq: i1} and \eqref{eq: i2}, we have
\begin{align*}
&\sqrt{k} \left(\widehat \Delta_n(s)-\widehat \Delta_n(s+1)-\delta(s)\right)\\
=&\sqrt{k} \left(\tilde\Delta_n\left(e_n,s\right) -\tilde\Delta_n\left(e_n,(s+1)\right)-\delta(s)\right)\\
=&\sqrt{k} \left(\tilde\Delta_n\left(e_n,s\right) -\Delta_n\left(e_n,s\right) \right)
-\sqrt{k} \left(\tilde\Delta_n\left(e_n,(s+1)\right) -\Delta_n\left(e_n,(s+1)\right) \right)\\
&+\sqrt{k} \left( \Delta_n\left(e_n,s\right) - \Delta_n\left(e_n,(s+1)\right)-\delta(s) \right)\\
=&  \nu_n(1, s)-\nu_n(1, s+1)-\delta(s) \nu_n(1,1) +o_p(1).
\end{align*}
Thus, to prove the result, it suffices to show that 
\begin{align}
\nu_n(1, s)-\nu_n(1, s+1)-\delta(s) \nu_n(1,1)\overset{d}{\rightarrow}&N(0, \sigma_2^2). \label{eq: vw1}
\end{align}
Define $I_i=\mathds{1}\left\{U_i<\frac{k}{n}< mU_{i+1,i+s-1}\right\}$, $J_i=\mathds{1}\left\{U_i<\frac{k}{n}< mU_{i+1,i+s}\right\}$ and 
$K_i=\mathds{1}\left\{U_i<\frac{k}{n}\right\}$. We need to show that
\begin{align*}
\frac{1}{\sqrt{k}}\sum_{i=1}^{n-s}\left(I_i-J_i-\delta(s) K_i-\E(I_i-J_i-\delta(s) K_i)\right) 
\overset{d}{\rightarrow} N(0, \sigma_2^2).
\end{align*}
We have
\begin{align}
&\var\left(\frac{1}{\sqrt{k}}\sum_{i=1}^{n-s}\left(I_i-J_i-\delta(s) K_i\right)\right)\nonumber\\
=&\frac{1}{k}\var\left(\sum_{i=1}^{n-s}I_i\right)+\frac{1}{k}\var\left(\sum_{i=1}^{n-s}J_i\right)
+\frac{1}{k}\delta^2(s)\var\left(\sum_{i=1}^{n-s}K_i\right)
-\frac{2}{k}\cov\left(\sum_{i=1}^{n-s}I_i, \sum_{i=1}^{n-s}J_i\right)\nonumber\\
&-\frac{2\delta(s)}{k}\cov\left(\sum_{i=1}^{n-s}I_i, \sum_{i=1}^{n-s}K_i\right)
+\frac{2\delta(s)}{k}\cov\left(\sum_{i=1}^{n-s}J_i, \sum_{i=1}^{n-s}K_i\right)\nonumber\\
=&\left(\Delta(s)+2\Lambda_2(1,1)\right)+\theta+\delta^2(s)(1+2\Lambda_1(1,1))+o(1) 
-\frac{2}{k}\cov\left(\sum_{i=1}^{n-s}I_i, \sum_{i=1}^{n-s}J_i\right)\nonumber\\
&-\frac{2\delta(s)}{k}\cov\left(\sum_{i=1}^{n-s}I_i, \sum_{i=1}^{n-s}K_i\right)
+2\delta(s)(\theta+\lambda_1),   \label{eq:var1}
\end{align}
where we used \eqref{eq: var_v1}, $\var(W_{s+1}(1))=\theta$, and $\frac{2}{k}\cov\left(\sum_{i=1}^{n-s}J_i, \sum_{i=1}^{n-s}K_i\right)=\theta+\lambda_1+o(1)$, from the proof for Theorem \ref{thm:AN}. Next, we compute two covariance term. By stationarity and Lemma \ref{lemma: cov_mixing}, 
\begin{align}
&\frac{1}{k}\cov\left(\sum_{i=1}^{n-s}I_i, \sum_{i=1}^{n-s}J_i\right)\nonumber\\
=&\frac{n}{k}\cov(I_1, J_1)+\frac{1}{k}\sum_{i=1}^{r_n} (n-s-i)\cov\left(I_1, J_{1+i}\right)
+\frac{1}{k}\sum_{i=1}^{r_n} (n-s-i)\cov\left(I_{1+i}, J_{1}\right)+o(1)\nonumber\\
=&\frac{n}{k}\E(I_1J_1)+\frac{n}{k}\sum_{i=1}^{r_n} \E\left(I_1 J_{1+i}\right)+\frac{n}{k}\sum_{i=1}^{r_n} E\left(I_{1+i}J_{1}\right)
+O\left(\frac{kr_n}{n}\right)+o(1)\nonumber\\
=&\frac{n}{k}\P\left(U_1<\frac{k}{n}< mU_{2,s+1}\right)
+\frac{n}{k}\sum_{i=1}^{r_n} \P\left(U_1<\frac{k}{n}< mU_{2,s},  U_{i+1}<\frac{k}{n}< mU_{i+2,i+s+1}\right)\nonumber\\
&+\frac{n}{k}\sum_{i=1}^{r_n} \P\left( U_{i+1}<\frac{k}{n}< mU_{i+2,i+s}, U_1<\frac{k}{n}< mU_{2,s+1}\right)
+o(1)\nonumber\\
=&\theta+\lambda_2+o(1),
 \label{eq:var2}
\end{align}
by Condition (A2) for $d=s+1$. Similarly, we have 
\begin{align}
&\frac{1}{k}\cov\left(\sum_{i=1}^{n-s}I_i, \sum_{i=1}^{n-s}k_i\right) \nonumber\\
=&\frac{n}{k}\P\left(U_1<\frac{k}{n}< mU_{2,s}\right)
+\frac{n}{k}\sum_{i=1}^{r_n} \P\left(U_1<\frac{k}{n}< mU_{2,s},  U_{i+1}<\frac{k}{n}\right)\nonumber\\
&+\frac{n}{k}\sum_{i=1}^{r_n} \P\left( U_{i+1}<\frac{k}{n}< mU_{i+2,i+s}, U_1<\frac{k}{n}\right)+o(1)\nonumber\\
=&\Delta(s)+\tilde \lambda_1+\lambda_3+o(1).   \nonumber\\
\label{eq:var3}
\end{align}
Combining \eqref{eq:var1}, \eqref{eq:var2} and \eqref{eq:var3}, it yields that 
\begin{align*}
&\var\left(\frac{1}{\sqrt{k}}\sum_{i=1}^{n-s}\left(I_i-J_i-\delta(s) K_i\right)\right)\\
\rightarrow
&\delta^2(s)\left(2\Lambda_1(1,1)-1\right)-2\delta(s)\left(\tilde\lambda_1-\lambda_1+\lambda_3-\frac{1}{2}\right)
+2\Lambda_2(1,1)-2\lambda_2=\sigma^2_2.
\end{align*}
And the Condition (2) in \cite{Utev1990} follows from the same argument as that for \eqref{eq: C2_Utev}. Therefore,
\eqref{eq: vw1} is proved.

\section{Lemmas and three propositions} \label{sec: lem}
We begin with two lemmas that are used for obtaining the covariance structures for $W_d$ in Proposition \ref{prop:AN}, $\tilde W$ in Proposition \ref{prop: marginal_AN} and $V(x)$ in the proof of Theorem \ref{thm:delta}. 
\begin{lemma}  \label{lemma: cov}
Define $ I_i(x):= \mathds{1}\left\{U_i<\frac{kx}{n}<  mU_{i+1,i+d-1}\right\} $ for $i=1,\ldots,n-d+1$ and $1/2\leq x\leq 3/2$.
Assume that $\frac{r_n k}{n}=o(1)$.  Let $y\in [1/2, 3/2]$.
\begin{itemize}
\item[(i)] If Condition A(2) holds, then,
\begin{equation*}
\sum_{i=1}^{r_n}(n-i)\cov(I_1 (x), I_{1+i}(x))=o(k), \label{eq: covII}
\end{equation*}
and 
\begin{equation*}
\sum_{i=d}^{r_n}(n-i)\cov(I_1(x), I_{1+i}(y))=o(k). \label{eq: covIJ}
\end{equation*}
\item[(ii)] If Condition (A3) holds, then 
\begin{align*}
\frac{1}{k}\sum_{i=1}^{r_n}(n-i)\cov(\mathds{1}\{U_1<{kx}/{n}\}, \mathds{1}\{U_{1+i}<{kx}/{n}\} )
\rightarrow\Lambda_1(x, x),
\end{align*}
and 
\begin{align*}
\frac{1}{k}\sum_{i=1}^{r_n}(n-i)\cov(\mathds{1}\{U_1<{kx}/{n}\}, \mathds{1}\{U_{1+i}<{ky}/{n}\} )
\rightarrow\Lambda_1(x, y).
\end{align*}
\item[(iii)] If the following limit exists:
\begin{align}
&\frac{n}{k}\sum_{i=d+1}^{r_n}\P(U_1<{kx}/{n}<mU_{2,d}, U_{i}<{ky}/{n}<mU_{i+1,i+d-1}) \nonumber\\
\rightarrow&\Lambda_2(x, y) \in [0,\infty), \label{eq: condition_a2d}
\end{align}
 then
\begin{align}
\frac{1}{k}\sum_{i=1}^{r_n}(n-i)\cov(I_1 (x), I_{1+i}(x))\rightarrow\Lambda_2(x, x), \nonumber\\
\label{eq: covII_2}
\end{align}
and 
\begin{align}
\frac{1}{k}\sum_{i=d}^{r_n}(n-i)\cov(I_1(x), I_{1+i}(y))=\Lambda_2(x, y).\nonumber\\
\label{eq: covIJ_2}
\end{align}
\end{itemize} 
\end{lemma}

\noindent{\bf Proof}   ~~Observes that results in parts (i) and (ii)  are both special cases of part (iii). For part (i), Condition A(2) implies that \eqref{eq: condition_a2d} holds and  $\Lambda_2(x, y)=0$, for $x, y \in [1/2, 3/2]$. Part(ii) is a special case of part(iii) with $d=1$. When $d=1$, $\Lambda_2$ in \eqref{eq: condition_a2d}
coincides with $\Lambda_1$ in Condition (A3).  Following we show the proof for part (iii).

Note that $\E(I_1(x))\leq \P(U_1<2k/n)=O\left(\frac{k}{n}\right)$, for any $x\in [1/2, 3/2]$. By construction, $\E(I_1(x)I_{1+i}(x))=0$ for $1\leq i\leq d-1$.
Thus, by \eqref{eq: condition_a2d},
\begin{align*}
&\frac{1}{k}\sum_{i=1}^{r_n}(n-i)\cov(I_1(x), I_{1+i}(x))
=\frac{1}{k}\sum_{i=1}^{r_n}(n-i)\left(\E(I_1(x)I_{1+i}(x))-(\E(I_1(x)))^2\right)\\
=&\frac{1}{k}\sum_{i=d}^{r_n}(n-i)\E(I_1(x)I_{1+i}(x))-\frac{(\E(I_1(x)))^2}{k}\sum_{i=1}^{r_n}(n-i)\\
\rightarrow& \Lambda_2(x, x).
\end{align*}
Hence \eqref{eq: covII_2} is proved. And \eqref{eq: covIJ_2} follows in the same way.
 \hfill  $\square$

\begin{lemma}\label{lemma: cov_mixing}
Let $A\in \sigma\left(\mathds{1}\left\{U_j \leq \frac{k}{n}\right\}, 1\leq j\leq d_1\right)$
and $B_i\in \sigma\left(\mathds{1}\left\{U_j \leq \frac{k}{n}\right\}, i\leq j\leq i+d_2\right)$, $i=1,2,\ldots,n-d_2$, where $d_1$ and $d_2$ are two positive integers. 
If $\frac{n}{k}\P(A)=O(1)$ and $\sum_{i=r_n}^n\left(1-\frac{i}{n}\right)\phi(i)=o(1)$, then 
\begin{equation}
\frac{1}{k}\sum_{i=r_n}^{n-d_2-1}(n-i)\cov(A,B_{i+1})=o(1).
\end{equation}
\end{lemma}
\noindent{\bf Proof}   ~~By the definition of $\phi_n(i)$ in  \eqref{eq: phi}, we have that  
$\lvert\cov(A,B_{i+1})\rvert=\lvert\P(A\cap B_{i+1})-\P(A)\P(B_{i+1})\rvert=\P(A)\lvert \P(B_{i+1}|A)-\P(B_{i+1})\rvert \leq \P(A)\phi_n(i+1-d_1)$. Thus,
\begin{align*}
&\frac{1}{k}\sum_{i=r_n}^{n-d_2-1}(n-i)\lvert\cov(A,B_{i+1})\rvert \\
\leq&\frac{n}{k}\P(A)\sum_{i=r_n}^{n-d_2-1} \left(1-\frac{i}{n}\right)\phi_n(i+1-d_i)=o(1). 
\end{align*}

\begin{proposition} \label{prop:d*}
Suppose that $d^u\geq d^*$ and  the assumptions of Corollary \ref{cor:delta} hold, then 
$$
\P\left(\hat d^*(k)=d^*\right)\overset{p}{\rightarrow} 1
$$
where $\hat d^*(k)$ is defined in \eqref{eq: hat_d}.
\end{proposition}
\noindent{\bf Proof}   ~~Denote $\widehat \Delta_n(i)=\widehat \Delta_n(i)-\widehat \Delta_n(i+1)$. By the definition of $\hat d^*(k)$,
$$
\{\hat d^*(k)\geq d^*+1\}\subset \left\{\max_{d^*\leq i\leq d^u}\widehat \Delta_n(i)\geq \frac{1}{\sqrt{k}}\right\},
$$
the latter has a probability tending to zero by Corollary \ref{cor:delta}.
On the other hand, for any $j\leq d^*-1$,
$$
\{\hat d^*(k)=j\}\subset \left\{\max_{j\leq i\leq d^u}\widehat \Delta_n(i)< \frac{1}{\sqrt{k}}\right\} \subset \left\{\widehat \Delta_n(d^*-1)< \frac{1}{\sqrt{k}}\right\},
$$
which also has  a probability tending to zero. Therefore,
$$
\P(\hat d^*(k)\neq d^*)\rightarrow 0.
$$

\begin{proposition} \label{prop: ARC}
For the AR(1) model with Cauchy margin defined in  \eqref{eq:ARC}, we have
\begin{itemize}
\item[(a)] for $z\geq 0$, $\ell_s(\vc_s)=s-(s-1)z$, for $s\geq 2$;
\item[(b)] for $z<0$, $\ell_2(\vc_2)=2$ and $\ell_s(\vc_s)= s-|z|^2$ for $s\geq 3$.
\end{itemize}
\end{proposition}
\noindent{\bf Proof}   ~~This result is easily derived by using the independence of $(X_1, \epsilon_2,\ldots, \epsilon_{s})$. 
Let $v_n$ be such that $\lim_{n\rightarrow \infty} n\P(\epsilon_i>v_n)=1$. 
Then $v_n=(1-|z|)u_n$.
 For $s\geq 2$, we have
\begin{align*}
	&\ell_s(\vc_s) = \lim_{n\to\infty}n\P\left(X_1>u_n \mbox{ or } \dots \mbox{ or } X_s>u_n\right)\\
	=&\lim_{n\to\infty}n \P(X_1>u_n \mbox{ or } \dots \mbox{ or } z^{s-1}X_1+z^{s-2}\epsilon_2
	+\cdots+\epsilon_{s}>u_n)\\
	=&\lim_{n\to\infty}n \P(\frac{X_1}{u_n}>1 \mbox{ or } \dots \mbox{ or } z^{s-1}\frac{X_1}{u_n}+z^{s-2}\frac{\epsilon_2}{u_n}
	+\cdots+\frac{\epsilon_{s}}{u_n}>1)\\
	=&\lim_{n\to\infty}n \P(\frac{X_1}{u_n}>1 \mbox{ or } \dots \mbox{ or } z^{s-1}\frac{X_1}{u_n}+z^{s-2}(1-|z|)\frac{\epsilon_2}{v_n}
	+\cdots+(1-|z|)\frac{\epsilon_{s}}{v_n}>1)\\
		=&\nu \{(t_1, \ldots,t_s): t_1>1 \mbox{ or }  \ldots \mbox{ or } z^{s-1}t_1+z^{s-2}(1-|z|)t_2
		+\cdots+(1-|z|)t_s>1\},
\end{align*}
 where $\nu$ denotes the exponent measure of $(X_1, \epsilon_2,\ldots, \epsilon_{s})$ (For definition of exponent measure, see Section 6.1.3 in \cite{deHaanFerreira2006}).  The last convergence follows from Theorem 6.1.11 in \cite{deHaanFerreira2006} and the fact that the distribution of $(X_1, \epsilon_2,\ldots, \epsilon_{s})$ belongs to the max domain of attraction. Now because of the exact independence of $X_1$ and the $\epsilon_i$'s and therefore asymptotic independence, the exponent measure $\nu$ puts mass only in the axes: $\nu\{(t_1, \ldots,t_s):t_i>a_1 \mbox{ and } t_j>a_2\}=0$ for any $i\neq j$ and positive $a_1$ and $a_2$. Then, the result readily follows from  the property that $\nu\{(t_1, \ldots,t_s):t_i>a_1 \}=1/a_1$.

\begin{proposition} \label{prop: ARCH}
An ARCH model defined in \eqref{eq:ARCH} does not satisfy $D^{(d)}(u_n)$ condition for any finite $d$. 
\end{proposition}
\noindent{\bf Proof}   ~~We apply Proposition 6.2 of \cite{Ehlertetal2015} to show that for any finite $d$,
$$
\lim_{n\rightarrow\infty}\P(M_{2,d}\leq u_n<M_{d+1,r_n}|X_1>u_n)>0.
$$ 
In this proof, all the cited equations  are referred to the formulas in \cite{Ehlertetal2015}. Note that ARCH(1,1) model is a special case of the model considered in that paper, which corresponds $\delta_1=\beta_1=0$ in the model given by relations (6.2) and (6.3) in that paper. Therefore $\phi(x)=\alpha_1^{-1/2}|x|$, for the $\phi$ appeared in the limit of (6.14) in that paper. Let $W$ denote a random variable from Pareto distribution with parameter $\alpha$ and $(Z_i)_{i\geq 1}$ are i.i.d. standard normal random variables. Then by Proposition 6.2 of \cite{Ehlertetal2015},
\begin{align*}
&\lim_{n\rightarrow\infty}\P(M_{2,d}\leq u_n<M_{d+1,r_n}|X_1>u_n)\\
\geq &\lim_{n\rightarrow\infty}\P(M_{2,d}\leq u_n<X_{d+1}|X_1>u_n)\\
=& \P\left(\max_{2\leq i\leq d}\frac{W}{|Z_0|}Z_i \prod_{j=0}^{i-1}\phi(Z_j)\leq 1 < \frac{W}{|Z_0|}Z_{d+1} \prod_{j=0}^{d}\phi(Z_j)\right)\\
=&\P\left(W \max_{2\leq i\leq d} \alpha_1^{i/2}Z_i \prod_{j=1}^{i-1}|Z_j|\leq 1 < W \alpha_1^{(d+1)/2}Z_{d+1} \prod_{j=1}^{d}|Z_j|\right)\\
\geq&\P\left(W> \alpha_1^{-(d+1)/2},\max_{2\leq i\leq d}Z_i\leq -1, Z_{d+1}>1 \right)\\
=&\alpha_1^{\alpha(d+1)/2}(\Phi(-1))^d>0,
\end{align*}
 where $\alpha_1 \in (0, 1)$ (which equals $1/2$ in our simulation example), $\alpha>0$ and $\Phi$ is  a standard normal distribution function.


\end{document}